\documentclass[11pt]{article}

\usepackage{amsmath}
\usepackage{amssymb,bbm}
\usepackage{mathabx}
 \usepackage{IEEEtrantools}
\usepackage{mathtools}

\usepackage{mathrsfs}
\usepackage{url}
\usepackage{hyperref}
\usepackage{xcolor}
\usepackage{color}

\usepackage{tikz-cd}

\usepackage{authblk}

\usepackage{multirow}

\usepackage{cases}

\usetikzlibrary{arrows,decorations.markings}
 \usetikzlibrary{fit,matrix}

\newcommand{\mychi}{\raisebox{.4ex}[1ex][0.5ex]{$\chi$}}

\newcommand{\prf}{\noindent{\it Proof\/}: }
\newcommand{\one}{\mathbbm{1}}
\newcommand{\id}{\text{id}}
\def \qed { \mbox{}\hfill
$\Box$\vspace{1ex}}

\newcommand{\uqslnc}{U_{q}(s\ell(n+1,\mathbb{C}))}
\newcommand{\tnsr}{\mbox{$\bigcirc\hspace{-0.80em}\mbox{\raisebox%
{-.43ex}{$\top$}}\;$}}

\newtheorem{dfn}{Definition}[section]
\newtheorem{thm}[dfn]{Theorem}
\newtheorem{theorem}[dfn]{Theorem}
\newtheorem{lmma}[dfn]{Lemma}
\newtheorem{lemma}[dfn]{Lemma}
\newtheorem{ppsn}[dfn]{Proposition}
\newtheorem{proposition}[dfn]{Proposition}
\newtheorem{crlre}[dfn]{Corollary}
\newtheorem{corollary}[dfn]{Corollary}
\newtheorem{xmpl}[dfn]{Example}

\newtheorem{rmrk}[dfn]{Remark}
\newtheorem{remark}[dfn]{Remark}
\newtheorem{xrcs}[dfn]{Exercise}

\newcommand{\bdfn}{\begin{dfn}\rm}
\newcommand{\bthm}{\begin{thm}}
\newcommand{\blmma}{\begin{lmma}}
\newcommand{\bppsn}{\begin{ppsn}}
\newcommand{\bcrlre}{\begin{crlre}}
\newcommand{\bxmpl}{\begin{xmpl}}
\newcommand{\brmrk}{\begin{rmrk}\rm}
\newcommand{\bxrcs}{\begin{xrcs}\rm\footnotesize}

\newcommand{\edfn}{\end{dfn}}
\newcommand{\ethm}{\end{thm}}
\newcommand{\elmma}{\end{lmma}}
\newcommand{\eppsn}{\end{ppsn}}
\newcommand{\ecrlre}{\end{crlre}}
\newcommand{\exmpl}{\end{xmpl}}
\newcommand{\ermrk}{\end{rmrk}}
\newcommand{\exrcs}{\end{xrcs}}

\numberwithin{dfn}{subsection}
\counterwithin*{equation}{section}

\begin{document}
\author[1]{\sc Manabendra Giri\thanks{manabgiri18r@isid.ac.in, manabendra991@gmail.com}}
\author[1]{\sc Arup Kumar Pal\thanks{arup@isid.ac.in, arupkpal@gmail.com}}
\affil[1]{Indian Statistical Institute, Delhi, INDIA}
\title{Quantized function algebras at $q=0$: type $A_{n}$ case}
\maketitle
 \begin{abstract} 
  We define the notion of  quantized function algebras  at $q=0$
  or crystallization of the $q$ deformations of the type $A_{n}$ compact Lie
  groups at the $C^*$-algebra level. The $C^{*}$-algebra $A_{n}(0)$ is defined
  as a universal $C^*$-algebra given by a finite set of generators and
  relations. We obtain these relations by looking at the irreducible
  representations of the quantized function algebras for $q>0$ and taking limit
  as $q\to 0+$ after rescaling the generating elements appropriately. We then
  prove that in the $n=2$ case the irreducible representations $A_{2}(0)$ are
  precisely the $q\to 0+$ limits of the irreducible representations of the
  $C^*$-algebras $A_{2}(q)$.
 \end{abstract}
{\bf AMS Subject Classification No.:}
17B37, 
20G42, 
46L67, 
\\





 \section{Introduction}
The theory of crystal bases was developed independently by  Kashiwara
(\cite{Kas-1990aa}, \cite{Kas-1991aa}, \cite{Kas-1991ab}, \cite{Kas-1993aa},
\cite{Kas-1994aa}, \cite{Kas-1995ab}) and Lusztig (\cite{Lus-1990aa},
\cite{Lus-1990ab}) in the early 1990's. For a broad class of Lie algebras
$\mathfrak{g}$, they observed that even though the quantized universal
enveloping algebras $U_q(\mathfrak{g})$ do not exist at $q=0$, when one looks at
a finite dimensional module over $U_q(\mathfrak{g})$, one can renormalize the
actions of the generating elements in such a way that they continue to make
sense at $q=0$. This was the starting point of the theory which later found
remarkable applications in representations of classical groups where one can use
the theory to construct good bases for finite dimensional representations of a
complex simple Lie algebra. This is considered to be one of the major
achievements of quantum group theory.

The coordinate ring of regular functions for the $q$-deformation of a complex
simple Lie group $G$ can be viewed as a module over the corresponding quantized
universal enveloping algebra $U_{q}(\mathfrak{g})$ of the Lie algebra
$\mathfrak{g}$ of the group $G$. Crystallization of these modules were studied
by Kashiwara in \cite{Kas-1993aa}. What we are interested in is the
crystallization of the $C^{*}$-algebraic counterpart of the coordinate rings,
i.e.\ the $C^{*}$-algebra of `continuous functions' for the $q$-deformation $G_
{q}$ of compact Lie groups $G$. These have not been investigated in the
literature. For the quantum $SU(2)$ group, Woronowicz had described the `$C^
{*}$-algebra $C(SU_{q}(2))$ at $q=0$' in \cite{Wor-1987aa}, which he obtained
by simply taking the relations that generate the $C^*$-algebra $C(SU_{q}
(2))$ and replacing $q$ by 0 in those relations and taking the universal
$C^*$-algebra generated by these new relations. This $C^*$-algebra did not seem
to appear again anywhere in the literature, except possibly in Hong \&
Szymanski (\cite{HonSzy-2002aa}) who observed that it is a graph $C^
{*}$-algebra.  However, very recently Chakraborty \& Pal  (\cite
{ChaPal-2022ey}) demonstrated a very interesting application of this $C^
{*}$-algebra in proving an approximate equivalence involving the GNS
representation of the Haar state of $SU_{q}(2)$ and the Soibelman
representation which is another well known faithful representation of the same
$C^{*}$-algebra. This approximate equivalence, in turn,  has interesting
connections with noncommutative  geometry and topology. In particular it gives
a $KK$-class via Cuntz's notion of quasihomomorphisms and also sets up a link
between the equivariant spectral triple for $SU_{q}(2)$ with another spectral
triple whose properties are easier to derive. In fact, it is this approximate
equivalence that is behind Connes' computations with the two spectral triples 
in \cite{Con-2004ab}. Thus the idea of looking at the $q\to 0$ limit 
can be traced back to the paper by Connes.
A similar equivalence
(or approximate equivalence) in the case of $SU_{q}(3)$ or more generally for
$q$-deformations of compact Lie groups should be an important tool in
understanding these quantum groups and their geometry better. Therefore it
seems worthwhile to investigate the crystallization of quantized function
algebras for the $q$-deformations of compact Lie groups in the $C^{*}$-algebra 
set up. In the present paper our goal is to initiate a study in
this direction. In particular, we will describe and study the $C^{*}$-algebras 
$C(G_{q})$ at $q=0$ where $G$ is a type $A_{n}$ compact Lie group.

For the rest of the paper, we will denote by $A_{n}(q)$ the $C^{*}$-algebra $C
(SU_{q}(n+1))$, where $q\in[0,1)$.  The $C^{*}$-algebras $A_{n}(q)$ for
$q\neq 0$ are given by a set of generators and relations. In order to obtain
the $C^{*}$-algebra $A_{n}(0)$, we obtain these relations `at
$q=0$' and then look at the universal $C^{*}$-algebra given by them. Our
first aim  is to obtain those relations. Recall that in Kashiwara's work,
while crystallizing an $U_{q}(\mathfrak{g})$ module, the passage to $q=0$ is
done through  quotienting by an appropriate ideal, but the actions of the
generating elements of $U_{q}(\mathfrak{g})$ need to be `renormalized' before
the quotienting procedure so that one gets something meaningful after
quotienting. We will use a limiting procedure instead of quotienting. Similar
to the Crystal bases case, we will make an appropriate scaling or
renormalization before taking limits. We do this in Section~2, where we look
at the ireducible representations of the $C^{*}$-algebra $A_{n}(q)$ for
$q\neq 0$ and examine the actions of the generating elements as $q\to
0+$. We arrive at a set of relations and prove the
existence of a universal $C^{*}$-algebra given by these relations. We call
this  universal $C^{*}$-algebra the crystallization of $A_{n}(q)$ and denote
it by $A_{n}(0)$. In Section~3, we focus on the $C^{*}$-algebra $A_{2}
(0)$ more closely. In particular, we  show that all the irreducible
representations arise exactly as limits of irreducible representations of $A_
{2}(q)$ as $q\to 0+$. As a simple application, we prove that $A_{2}(0)$ has a
natural coproduct making it a compact quantum semigroup, but it is not a
compact quantum group.

After the first version of this paper came out on the arXiv, Matassa \& Yuncken
(\cite{MatYun-2022aa}) introduced crystallizations of the $C^*$-algebras $C(G_
{q})$ for a broader class of compact Lie groups using a slightly different 
technique. This will be briefly discussed at the
end of the next section.

  \section{Crystallization of $A_{n}(q)$}

  \subsection{The quantum group $SU_{q}(n+1)$}
We will start by recalling a few facts on the quantum group $SU_{q}
(n+1)$. Associated to this quantum group, there are two main objects --
 the quantized universal enveloping algebra (QUEA) 
$U_{q}(sl(n+1,\mathbb{C}))$ (together with a $*$-structure called its 
compact real form)
and the quantized function algebra (QFA) $\mathcal{O}_{q}(SU(n+1))$. 
Our main object of interest will be the enveloping $C^*$-algebra 
$C(SU_{q}(n+1))$ of $\mathcal{O}_{q}(SU(n+1))$.
We recall a few facts about these  below.
For more details, we refer the reader to (\cite{Koe-1991aa}),
(Chapters~3--5, \cite{ParWan-1991kr}), (Chapters~8--9, \cite
{KliSch-1997aa}) and (Chapter 3, \cite{KorSoi-1998ab}).

Let $q\in\mathbb{R}_{+}$, $q\not\in\{0,1\}$. The QUEA $U_{q}(sl(n+1,\mathbb{C}))$
is a complex algebra generated by $K_i$, $K_i^{-1}$, $X_i$ and $Y_i$ for 
$i=1,\ldots,n+1$ obeying the following relations:
\begin{IEEEeqnarray}{rClClrCrCl}
   K_i K_i^{-1} &=&K_i^{-1} K_i &=& 1,&&\qquad \qquad&
   K_i K_j &=& K_j K_i   \label{eq:uqrel-1}\\ 
   K_i X_j K_i^{-1} &=& \IEEEeqnarraymulticol{3}{l}{q^{a_{i j}} X_j }&&&
   K_i Y_j K_i^{-1} &=& q^{-a_{i j}} Y_j \label{eq:uqrel-2}\\
   X_i Y_j-Y_j X_i &=& \IEEEeqnarraymulticol{8}{l}{\delta_{i j} \frac{K_i-K_i^{-1}}{q^{}-q^{-1}}}\label{eq:uqrel-3}\\
   \IEEEeqnarraymulticol{8}{r}{
   \sum_{k=0}^{1-a_{i j}}(-1)^k\left[\begin{array}{c}
1-a_{i j} \\
k
\end{array}\right]_{q} X_i^{1-a_{i j}-k} X_j X_i^k 
}
&=& 0  \quad (i \neq j) \label{eq:uqrel-4}\\
\IEEEeqnarraymulticol{8}{r}{
\sum_{k=0}^{1-a_{i j}}(-1)^k\left[\begin{array}{c}
1-a_{i j} \\
k
\end{array}\right]_{q} Y_i^{1-a_{i j}-k} Y_j Y_i^k 
}&=& 0 \quad (i \neq j)\label{eq:uqrel-5}
\end{IEEEeqnarray}
where $(\!(a_{i,j})\!)$ is the Cartan matrix for type $A_{n}$ simple Lie groups.
The following $*$-structure on $U_{q}(sl(n+1,\mathbb{C}))$ is called its compact real form:
\[
K_{i}^{*} = K_{i}^{}, \qquad
X_{i}^{*} = K_{i}^{}Y_{i}^{},\qquad
Y_{i}^{*} = X_{i}^{}K_{i}^{-1}.
\]

The space of matrix elements $\pi_{f,v}$ of all finite dimensional admissible (or type I)
representations $\pi$ (on $V^{\pi}$) of $\uqslnc$ inherits a Hopf algebra structure from
$\uqslnc$ via the pairing
$\langle \pi_{f,v}, a\rangle := f(\pi(a)v)$ where $v\in V^{\pi}$ and $f$ belongs to the dual of $V^{\pi}$.
This Hopf algebra is called the quantized algebra of regular functions on 
$SL(n+1,\mathbb{C})$ and will be denoted by $\mathcal{O}_{q}(SL(n+1,\mathbb{C}))$..

Let $u^{\one}$ denote the fundamental representation of $\uqslnc$ corresponding to the 
weight $\one:=(1,0,\ldots,0)$. 
The $R$-matrix for this representation is given by
\[
 R_{q}
    =q^{-1}\sum_{i} E_{i,i}\otimes E_{i,i} 
    + \sum_{i\neq j} E_{i,i}\otimes E_{j,j} 
      + (q^{-1}-q)\sum_{i>j} E_{i,j}\otimes E_{j,i}.
\]
Let
$\sigma$ be the flip operator. Then $\hat{R}_{q}:=\sigma\circ R_{q}$
intertwines $u^{\one}\tnsr u^{\one}$ with itself, where
$u^{\one}\tnsr u^{\one}$ denotes the tensor product of the representation
$ u^{\one}$ with itself. Fix a highest weight vector $e_{1}$ and let
$\{e_{1},e_{2},\ldots,e_{n+1}\}$ be the basis given by 
$e_{j+1}=u^{\one}(Y_{j})e_{j}$.
Let $(\!(u_{i,j})\!)$ denote the matrix of $u^{\one}(\cdot)$ 
with respect to this basis
and its dual basis. The intertwiner equation
$\hat{R}_{q}( u^{\one}\tnsr u^{\one})(\cdot)=(u^{\one}\tnsr u^{\one})(
\cdot)\hat{R}_{q}$ gives the following commutation relations
(Equations~(2.1--2.2), \cite{Koe-1991aa}):
\begin{IEEEeqnarray}{rCll}
	 u_{i,j}u_{i,l} &=& qu_{i,l}u_{i,j} &
	  \text{if }j<l,\label{eq:comm-t1}\\
	 u_{i,j}u_{k,j} &=& q u_{k,j}u_{i,j} &
	  \text{if }i<k,\label{eq:comm-t2}\\
	 u_{i,l}u_{k,j} &=& u_{k,j}u_{i,l}
   & \text{if }i<k\text{ and }j<l,\label{eq:comm-t3}\\
	 u_{i,j}u_{k,l} - u_{k,l}u_{i,j} &=& (q-q^{-1})u_{i,l}u_{k,j}
	  \quad &\text{if }i<k\text{ and }j<l.\label{eq:comm-t4}
\end{IEEEeqnarray}
The quantum determinant for the matrix $(\!(u_{i,j})\!)$ is given by (see
Equation~(2.7), \cite{Koe-1991aa} or Chapter~4, \cite{ParWan-1991kr})
\begin{equation}\label{eq:comm-t5}
	D_{q}=\sum_{\sigma\in\mathscr{S}_{n+1}}(-q)^{\ell(\sigma)}
      u_{1,\sigma(1)}u_{2,\sigma(2)}\dots u_{n+1,\sigma(n+1)},
\end{equation}
where $\ell(\sigma)$ denotes the length of the permutation $\sigma$.
Let $D_{q}^{s,r}$ denote the $(s,r)$\raisebox{.2ex}{th} cofactor of the matrix
$(\!(u_{i,j})\!)$. Then one has
\[
(\!(u_{i,j})\!)(\!((-q)^{r-s}D_{q}^{s,r})\!)=
(\!((-q)^{r-s}D_{q}^{s,r})\!)(\!(u_{i,j})\!)=D_{q}I.
\]
It is also known that the quantum determinant $D_{q}$ is a central element of
the algebra generated by the relations (\ref{eq:comm-t1}-\ref{eq:comm-t4})
(Theorem~4.6.1, page-50, \cite{ParWan-1991kr}). It follows from this
centrality result that
\[
D_{q}^{s,r}u_{i,j}^{}=u_{i,j}^{}D_{q}^{s,r}\qquad\text{if } i\neq r,\; j\neq s.
\]
It is known that the quantized algebra $\mathcal{O}_{q}(SL(n+1,\mathbb{C}))$
is the algebra generated by the relations (\ref{eq:comm-t1}-\ref{eq:comm-t4})
together with the relation
\begin{align}
  D_{q} &=1.\label{eq:qdet=1}
\end{align}
The compact real form of $\uqslnc$ induces a $*$-structure on  
$\mathcal{O}_{q}(SL(n+1,\mathbb{C}))$ that makes the matrix $(\!(u_{i,j})\!)$ 
unitary and therefore one has
\begin{align}
	u_{r,s}^{*} &= (-q)^{s-r}D_{q}^{r,s}.\label{eq:comm-t6}
\end{align}
Consequently, one also has
\begin{equation}\label{eq:comm-t7}
	u_{i,j}^{}u_{r,s}^{*}=u_{r,s}^{*}u_{i,j}^{}, \qquad i\neq r,\; j\neq s.
\end{equation}
The algebra $\mathcal{O}_{q}(SL(n+1,\mathbb{C}))$ together with this $*$-structure
is the Hopf $*$-algebra of coordinate functions on the quantum $SU(n+1)$ group.
We will denote it by $\mathcal{O}_{q}(SU(n+1))$. Its enveloping $C^*$-algebra 
is the $C^*$-algebra $A_{n}(q)\equiv C(SU_{q}(n+1))$ associated with $SU_{q}(n+1)$.

Next let us briefly recall the irreducible representations of the
$C^{*}$-algebra $A_{n}(q)$. The Weyl group for $SU_q(n+1)$ is isomorphic
to the permutation group $\mathscr{S}_{n+1}$ on $n+1$ symbols. Denote by
$s_i$ the transposition $(i,i+1)$. Then $\{s_{1},s_{2},\ldots,s_{n}\}$ form a
set of generators for $\mathscr{S}_{n+1}$. For $a,b\in\mathbb{N}$, $a\leq
b$, denote by $s_{[a,b]}$ the product $s_{b}s_{b-1}\ldots s_{a}$. Let
$\omega=s_{[a_{1},b_{1}]}s_{[a_{2},b_{2}]}\ldots s_{[a_{k},b_{k}]}$, where
$1\leq b_{1}<b_{2}<\ldots < b_{k}\leq n$. Then $\omega$ is a reduced word
in $\mathscr{S}_{n+1}$.

Let $S$  be the left shift operator and $N$ be the 
number operator on $\ell^{2}(\mathbb{N})$:
\[
Se_n=\begin{cases}
       e_{n-1}& \text{if } n\geq 1,\cr
       0 & \text{if }n=0.
	  \end{cases}, \qquad Ne_{k}=ke_{k},\quad k\in\mathbb{N}.
\]
Denote by $\psi_{s_k}^{(q)}$ the following representation of $A_{n}(q)$ on
$\ell^{2}(\mathbb{N})$:
\begin{equation}\label{eq:irr-rep-1}
\psi_{s_k}^{(q)}(u_{i,j}^{})=\begin{cases}
                   S\sqrt{I-q^{2N}} &\text{if }i=j=k,\cr
                  \sqrt{I-q^{2N}}S^{*} & \text{if } i=j=k+1,\cr
                   -q^{N+1} &  \text{if } i=k, j=k+1,\cr
                    q^{N}   &  \text{if }i=k+1, j=k,\cr
                    \delta_{i,j}I & \text{otherwise}.
					  \end{cases}
\end{equation}
For a reduced word $\omega=s_{i_1}s_{i_2}\ldots
s_{i_k}\in\mathscr{S}_{n+1}$, define $\psi_{\omega}^{(q)}$ to be
$\psi_{s_{i_1}}^{(q)}\ast\psi_{s_{i_2}}^{(q)}\ast\ldots\ast\psi_{s_{i_k}}^{
(q)}$. Here, for two representations $\phi$ and $\psi$, $\phi\ast\psi$
denote the representation $(\phi\otimes\psi)\Delta_{q}$, where $\Delta_{q}$ denotes the comultiplication map from $A_{n}(q)$ to $A_{n}(q)\otimes A_{n}(q)$.

Next, let $\lambda\equiv(\lambda_{1},\ldots,\lambda_{n})\in (S^1)^{n}$.
Define
\begin{equation}\label{eq:irr-rep-2}
	\mychi_{\lambda}(u_{i,j}^{})=
	\begin{cases}
            \lambda_{i}\delta_{i,j} &  \text{if }i=1,\cr
            \bar{\lambda}_{n} \delta_{i,j} &  \text{if }i=n+1,\cr
            \bar{\lambda}_{i-1}\lambda_{i}\delta_{i,j} &  \text{otherwise}.
	\end{cases}
\end{equation}
It is a well-known result of Soibelman 
(Theorem 6.2.7, page~121, \cite{KorSoi-1998ab}) that for any
reduced word $\omega\in\mathscr{S}_{n+1}$ and a tuple $\lambda\equiv
(\lambda_{1},\ldots,\lambda_{n})\in (S^1)^{n}$, the representation
\begin{equation}\label{eq:Aq(n)-irred}
  \psi_{\lambda,\omega}^{(q)}:=\mychi_{\lambda}\ast \psi_{\omega}^{(q)}
\end{equation}
is an irreducible representation of $A_{n}(q)$, and these give all the
irreducible representations of $A_{n}(q)$. Note that $\psi_{\lambda,\omega}^{(q)}$
acts on the Hilbert space $\ell^{2}(\mathbb{N}^{k})$,
where $k$ is the length of the reduced word $\omega=s_{i_1}s_{i_2}\ldots
s_{i_k}$ in $\mathscr{S}_{n+1}$.

  \subsection{The limiting relations}
In the remaining sections, we will denote the generating elements 
$u_{i,j}$ by $u_{i,j}(q)$ as we will be dealing with the algebras
$A_{n}(q)$ for different values of $q$.
Observe that for all $q\in (0,1)$, the irreducible representations 
$\psi_{\lambda,\omega}^{(q)}$ of the
$C^{*}$-algebras $A_{n}(q)$ are parametrized by the same set 
$(S^1)^{n}\times \mathscr{S}_{n+1}$ which is independent of $q$, 
and the representation $\psi_{\lambda,\omega}^{(q)}$ acts on the Hilbert space
$\mathcal{H}_{\lambda,\omega}:= \ell^{2}(\mathbb{N}^{k})$, $k$ being the length
of $\omega$, which does not change for
different values of the parameter $q$.
Thus for any irreducible representation $\psi_{\lambda,\omega}^{(q)}$ of $A_{n}(q)$,
the operators $\psi_{\lambda,\omega}^{(q)}(u_{i,j}(q))$ live on the same Hilbert space
$\mathcal{H}_{\lambda,\omega}$ for $q\in (0,1)$.
The following is the key observation in obtaining the  commutation 
relations used to define the crystallized $C^{*}$-algebra $A_{n}(0)$.
\begin{proposition}\label{prop:key-prop}
For any irreducible representation $\pi$ of $A_{n}(q)$ on a Hilbert space
$\mathcal{H}$, the norm limits $\lim_{q\to 0+}\pi(u_{i,j}(q))$ exist for all
$1\leq i,j\leq n+1$.

Moreover, for $i<j$, the norm limits
$\lim_{q\to 0+}(-q)^{i-j}\pi(u_{i,j}(q))$ exist.
\end{proposition}
\prf
Observe that if the above limits exist for two representations $\pi_{1}$
and $\pi_{2}$, then the limits also exist for the representation
$\pi_{1}\ast\pi_{2}$. This is immediate from the equality
\[
\pi_{1}\ast\pi_{2}(u_{i,j}^{}(q))=\sum_{k=1}^{n+1}\pi_{1}(u_{i,k}^{}(q))
\otimes \pi_{2}(u_{k,j}^{}(q)).
\]
Therefore the proof reduces to a simple verification of the existence of
the limits for the representations $\mychi_{\lambda}$ and $\psi_{s_{k}}$
using (\ref{eq:irr-rep-1}) and (\ref{eq:irr-rep-2}).
\qed

\begin{remark}
\begin{enumerate}
  \item
  From the second part of the result, it follows that for $i<j$, one has
  $\lim_{q\to 0+}\pi(u_{i,j}(q))=0$.
  \item
  The two statements in the above result can be combined into the following:
\begin{quote}
  the limit $\lim_{q\to 0+}(-q)^{\min\{i-j,0\}}\pi(u_{i,j}(q))$ exists and is
  finite for each ireducible representation $\pi$ of $A_{n}(q)$.
\end{quote}
\end{enumerate}
\end{remark}

 Notice that for different values of the parameter $q$, the elements $u_{i,j}
 (q)$ belong to different algebras. Therefore even though we can talk about
 the  limits 
 \[
 \lim_{q\to 0+}(-q)^{\min\{i-j,0\}}\pi(u_{i,j} (q))
 \]
 in the above proposition, it is not immediate how to
 interpret them conceptually. In order to do so, 
 one needs to bring in the rational form of
 the quantized function algebra $\mathcal{O}_{q}(SU(n+1))$.
  The rational form of the QUEA is an algebra $U_{t}(sl(n+1,\mathbb{C}))$ over
  $\mathbb{Q}(t)$ given by the relations (\ref{eq:uqrel-1}--\ref
  {eq:uqrel-5}) with $q$ replaced by the indeterminate $t$. The representation
  theory of $U_{t}(sl(n+1,\mathbb{C}))$ is parallel to that of $\uqslnc$, and
  the rational form $\mathcal{O}_{t}(SL(n+1,\mathbb{C}))$ of the QFA 
  is  the space of matrix elements
  $\pi_{v,f}$ of all finite dimensional admissible (or type I) representations
  $\pi$ of $U_{t}(sl(n+1,\mathbb{C}))$ on $V^{\pi}_{t}$ 
  (Proposition 7.2.2, \cite{Kas-1993aa}).

Let us denote the subring $\mathbb{Q}[t,t^{-1}]$ of $\mathbb{Q}(t)$ by $A$.
Let $U_{t}^{A}(sl(n+1,\mathbb{C}))$ denote the $A$-subring of 
$U_{t}(sl(n+1,\mathbb{C}))$ generated by the elements
\[
X_{i}^{(r)}:=\frac{X_{i}^{r}}{[r]_{t}!},\quad
   Y_{i}^{(r)}:=\frac{Y_{i}^{r}}{[r]_{t}!},\quad
   K_{i}^{\pm 1},
\]
where $[r]_{t}!$ stands for the $t$-factorial $\prod_{j=1}^{r}[j]_{t}$, with 
$[j]_{t}=\frac{t^{j}-t^{-j}}{t-t^{-1}}$. One then has (see Section 4.1, \cite{Kas-1993aa})
\[
U_{t}(sl(n+1,\mathbb{C}))\cong \mathbb{Q}(t)\otimes_{A}U_{t}^{A}(sl(n+1,\mathbb{C})).
\]
Corresponding to $U_{t}^{A}(sl(n+1,\mathbb{C}))$, one 
now defines an $A$-subring of $\mathcal{O}_{t}(SL(n+1,\mathbb{C}))$ as follows
(see equation~7.3.5, \cite{Kas-1993aa}):
\[
\mathcal{O}_{t}^{A}(SL(n+1,\mathbb{C})):=
  \bigl\{a\in \mathcal{O}_{t}(SL(n+1,\mathbb{C})):\langle a,x\rangle\in A
     \text{ for all }x\in   U_{t}^{A}(sl(n+1,\mathbb{C})) \bigr\}.
\]
One has the isomorphism (see Lemma~7.3.1, \cite{Kas-1993aa})
\[
\mathcal{O}_{t}(SL(n+1,\mathbb{C}))\cong 
     \mathbb{Q}(t)\otimes_{A} \mathcal{O}_{t}^{A}(SL(n+1,\mathbb{C})).
\]
For $q\in (0,\infty)$, 
one then has the specialization
\begin{equation}\label{eq:key-iso}
\mathcal{O}_{q}(SL(n+1,\mathbb{C}))\cong 
     \mathbb{C}\otimes_{A} \mathcal{O}_{t}^{A}(SL(n+1,\mathbb{C})).
\end{equation}
via the action of $A$ on $\mathbb{C}$ by 
$p[t,t^{-1}]\cdot\lambda=p[q,q^{-1}]\lambda$.

Note that $U_{t}^{A}(sl(n+1,\mathbb{C}))$ is closed under the $*$-structure
on $U_{t}(sl(n+1,\mathbb{C}))$ given by
\[
K_{i}^{*} =K_{i}^{}, \qquad
X_{i}^{*} = K_{i}^{}Y_{i}^{},\qquad
Y_{i}^{*} = X_{i}^{}K_{i}^{-1},\qquad
t^{*}=t.
\]
This induces a $*$-structure on $\mathcal{O}_{t}(SL(n+1,\mathbb{C}))$.
We will denote the resulting $*$-algebra over $\mathbb{Q}(t)$ by $\mathcal{O}_{t}(SU(n+1))$.
The $A$-subring $\mathcal{O}_{t}^{A}(SL(n+1,\mathbb{C}))$
is closed under the $*$-structure on $\mathcal{O}_{t}(SL(n+1,\mathbb{C}))$.
We will denote by 
$\mathcal{O}_{t}^{A}(SU(n+1))$ the $A$-subring 
$\mathcal{O}_{t}^{A}(SL(n+1,\mathbb{C}))$ together with 
this inherited $*$-structure. 

The isomorphism (\ref{eq:key-iso}) now allows one to get a map 
\[
\theta_{q}:\mathcal{O}_{t}^{A}(SU(n+1))\to \mathcal{O}_{q}(SU(n+1))\subseteq A_{n}(q)
\]
that sends $p[t,t^{-1}]u_{i,j}(t)$ to $p[q,q^{-1}]u_{i,j}(q)$ for each $i$, $j$ and for 
any polynomial $p[t,t^{-1}]$.
In particular, $(-t)^{\min\{i-j,0\}}u_{i,j}(t)\in \mathcal{O}_{t}^{A}(SU(n+1))$ and
one has 
\[
(-q)^{\min\{i-j,0\}}u_{i,j}(q)=\theta_{q}\bigl((-t)^{\min\{i-j,0\}}u_{i,j}(t)\bigr)
\]
for each $q\in (0,1)$.
Thus the above proposition says that for each irreducible 
representation $\psi_{\lambda,\omega}^{(q)}$ of $A_{n}(q)$, the norm limits
$\lim_{q\to 0+}\psi_{\lambda,\omega}^{(q)}\circ\theta_{q}
         \bigl((-t)^{\min\{i-j,0\}}u_{i,j}(t)\bigr)$ 
exist and are finite.


We now proceed to find the relations satisfied by the limiting operators.
For an irreducible representation $\pi$ of $A_{n}(q)$, let us define
\begin{align}
Z_{i,j}^{\pi} &:=
  \begin{cases}
	  \lim_{q\to 0+}\pi(u_{i,j}^{}(q)) & \text{if } i\geq j,\\
	  \lim_{q\to 0+}(-q)^{i-j}\pi(u_{i,j}^{}(q)) & \text{if } i < j.
  \end{cases}\label{eq:def-y}
\end{align}

\begin{proposition}
For any irreducible representation $\pi$ of $A_{n}(q)$, the operators
$Z_{i,j}\equiv Z_{i,j}^{\pi}$ satisfy the following commutation relations:
\begin{IEEEeqnarray}{rClll}
	 Z_{i,j}^{}Z_{i,l}^{} &=& 0 & \text{if}\quad
	          &j<l,\label{eq:comm-y1}\\
	 Z_{i,j}^{}Z_{k,j}^{} &=& 0 & \text{if}\quad
	          &i<k,\label{eq:comm-y2}\\
            &&&&\nonumber\\
     Z_{i,l}^{}Z_{k,j}^{}  &=&0 & \text{if}\quad
	          &\begin{cases}
		      j<l\leq i<k,\\
	            i<k\leq j<l.
			\end{cases}\label{eq:comm-y3}\\
      &&&&\nonumber\\
	  Z_{i,l}^{}Z_{k,j}^{} - Z_{k,j}^{}Z_{i,l}^{} &=& 0 \quad
	   & \text{if}\quad &i<k\text{ and }j<l.\label{eq:comm-y4}\\
     &&&&\nonumber\\
  Z_{i,j}^{}Z_{k,l}^{} - Z_{k,l}^{}Z_{i,j}^{}  &=& 0 &\text{if}\quad
	  &\begin{cases}
	          i<j<j+1<k<l,\\
			  i<j<j+1<l\leq k,\\
			    j\leq i<i+1<k<l,\\
			    j\leq i<i+1<l\leq k,
	  \end{cases}\label{eq:comm-y5}\\
    &&&&\nonumber\\
Z_{i,j}^{}Z_{k,l}^{} - Z_{k,l}^{}Z_{i,j}^{}&=& Z_{i,l}^{}Z_{k,j}^{}\quad
     & \text{if}\quad &
	 \begin{cases}
	   i<j<j+1=k<l,\\
       i<j<j+1=l \leq k,\\
       j\leq i<i+1=k<l,\\
       j\leq i<i+1=l\leq k.
	  \end{cases}\label{eq:comm-y6}
\end{IEEEeqnarray}
\end{proposition}
\prf
The proof is a careful verification in each of the cases using the
relations (\ref{eq:comm-t1}--\ref{eq:comm-t4}) and
Proposition~\ref{prop:key-prop}. First observe that the $\pi(u_{i,j}^{})$'s
satisfy the same relations (\ref{eq:comm-t1}--\ref{eq:comm-t4}) as the
$u_{i,j}^{}$'s. Therefore the relations (\ref{eq:comm-y1}),
(\ref{eq:comm-y2}) and (\ref{eq:comm-y4}) follow from
(\ref{eq:comm-t1}--\ref{eq:comm-t3}). Similarly the relations in
(\ref{eq:comm-y3}) follow from (\ref{eq:comm-t4}).
Let us next prove the first relation in (\ref{eq:comm-y5}). Let $i<j<j+1<k<l$.
Note that in this case we have
\begin{IEEEeqnarray*}{rClrCl}
	Z_{i,j}^{} &=& \lim_{q\to 0+}(-q)^{i-j}\pi(u_{i,j}^{}(q)),\qquad &
	  Z_{k,l}^{} &=& \lim_{q\to 0+}(-q)^{k-l}\pi(u_{k,l}^{}(q)),\\
	Z_{i,l}^{} &=& \lim_{q\to 0+}(-q)^{i-l}\pi(u_{i,l}^{}(q)),\quad &
	  Z_{k,j}^{} &=& \lim_{q\to 0+}\pi(u_{k,j}^{}(q)).
\end{IEEEeqnarray*}
From (\ref{eq:comm-t4}), we have
\begin{IEEEeqnarray*}{rCl}
	\IEEEeqnarraymulticol{3}{l}{
	(-q)^{i-j}u_{i,j}^{}(q)(-q)^{k-l}u_{k,l}^{}(q) - (-q)^{k-l}u_{k,l}^{}(q)(-q)^{i-j}u_{i,j}^{}(q)}\\
	\qquad \qquad\qquad&=& (-q)^{i-j+k-l}(q-q^{-1})u_{i,l}^{}(q)u_{k,j}^{}(q)\\
	&=& (-q)^{k-j}(q-q^{-1})q^{i-l}u_{i,l}^{}(q)u_{k,j}^{}(q)\\
	&=& (-1)^{k-j}(q^{k-j+1}-q^{k-j-1})q^{i-l}u_{i,l}^{}(q)u_{k,j}^{}(q).
\end{IEEEeqnarray*}
Since $k-j-1>0$, it follows that
$Z_{i,j}^{}Z_{k,l}^{} - Z_{k,l}^{}Z_{i,j}^{}  = 0$.
Proofs of the remaining relations in (\ref{eq:comm-y5}) are similar.

Next we will prove the first relation in (\ref{eq:comm-y6}).
Let $i<j<j+1=k<l$. From (\ref{eq:comm-t4}), as before we have
\begin{IEEEeqnarray*}{rCl}
	\IEEEeqnarraymulticol{3}{l}{
  (-q)^{i-j}u_{i,j}^{}(q)(-q)^{k-l}u_{k,l}^{}(q) - (-q)^{k-l}u_{k,l}^{}(q)(-q)^{i-j}u_{i,j}^{}(q)}\\
	\qquad \qquad\qquad \qquad \qquad\qquad&=& (-1)^{k-j}(q^{k-j+1}-q^{k-j-1})q^{i-l}u_{i,l}^{}(q)u_{k,j}^{}(q).
\end{IEEEeqnarray*}
In this case, since $k-j=1$, we have
 $Z_{i,j}^{}Z_{k,l}^{} - Z_{k,l}^{}Z_{i,j}^{}=
 	       Z_{i,l}^{}Z_{k,j}^{}$. Proofs of the remaining relations are similar.
\qed

\begin{remark}
The relations (\ref{eq:comm-y1}--\ref{eq:comm-y6}) above can be rewritten in a slightly more compact form as follows
\begin{IEEEeqnarray}{rClll}
   Z_{i,j}^{}Z_{i,l}^{} &=& 0 & \text{if }
            &j<l,\label{eq:comm-w1}\\
   Z_{i,j}^{}Z_{k,j}^{} &=& 0 & \text{if }
            &i<k,\label{eq:comm-w2}\\
            &&&&\nonumber\\
    Z_{i,l}^{}Z_{k,j}^{} - Z_{k,j}^{}Z_{i,l}^{} &=& 0 \quad
     & \text{if } &i<k\text{ and }j<l.\label{eq:comm-w4}\\
     &&&&\nonumber\\
     Z_{i,l}^{}Z_{k,j}^{}  &=&0 & \text{if }& i<k,\,j<l \text{ and }
           \max\{i,j\}\geq \min\{k,l\},\label{eq:comm-w3}\\
      &&&&\nonumber\\
  Z_{i,j}^{}Z_{k,l}^{} - Z_{k,l}^{}Z_{i,j}^{}  &=& 0 &\text{if }&i<k,\,j<l \text{ and }
           \max\{i,j\}+1 < \min\{k,l\},\label{eq:comm-w5}\\
    &&&&\nonumber\\
Z_{i,j}^{}Z_{k,l}^{} - Z_{k,l}^{}Z_{i,j}^{}&=& Z_{i,l}^{}Z_{k,j}^{}\quad
     & \text{if }&i<k,\,j<l \text{ and }
           \max\{i,j\}+1= \min\{k,l\},\label{eq:comm-w6}
\end{IEEEeqnarray}
\end{remark}

\begin{proposition}\label{ppsn:prod=1}
For any irreducible representation $\pi$ of $A_{n}(q)$, one has
\begin{align}
  Z_{1,1}^{} Z_{2,2}^{}\dots  Z_{n+1,n+1}^{} &= I.\label{eq:prod=1}
\end{align}
\end{proposition}
\prf
From (\ref{eq:qdet=1}) one gets $\pi(D_{q})=I$ and the result follows by taking limit as $q\to 0+$. 
\qed

From (\ref{eq:comm-t6}), it follows that for any representation $\pi$ of
$A_{n}(q)$, one has
\begin{align}
(\pi(u_{r,s}(q)))^{*} &= (-q)^{s-r}\pi(D_{q}^{r,s}).\label{eq:comm-ystar}
\end{align}
We will need the following lemma in the proof of the next proposition.
\begin{lemma} \label{lmma:perm}
Let $1\leq s < r \leq n+1$. Let
\[
i_{k}=\begin{cases}
          k & \text{if } 1\leq k < r,\\
		  k+1 &\text{if } r\leq k \leq n,
	  \end{cases} \qquad
j_{k}=\begin{cases}
          k & \text{if } 1\leq k < s,\\
		  k+1 &\text{if } s\leq k \leq n,
	  \end{cases}
\]
Let $\sigma\in\mathscr{S}_{n}$. Then
\[
\sum_{k: j_{\sigma(k)}>i_{k}}(j_{\sigma(k)}-i_{k})\geq r-s.
\]
\end{lemma}
\prf
The required inequality follows from the following simple computation:
\begin{align*}
\sum_{k: j_{\sigma(k)}>i_{k}}(j_{\sigma(k)}-i_{k})
   &\geq \sum_{k=1}^{n}(j_{\sigma(k)}-i_{k})\\
   &= \sum_{k=1}^{n}j_{\sigma(k)} -\sum_{k=1}^{n}i_{k}\\
   &= \left(\frac{(n+1)(n+2)}{2}-s\right)-\left(\frac{(n+1)(n+2)}{2}-r\right)\\
   &= r-s.
\end{align*}
\qed

\begin{proposition}
For any irreducible representation $\pi$ of $A_{n}(q)$, one has
\begin{align}
	Z_{r,s}^{*}
	  &= \begin{cases}
	       Z_{1,1}^{}\dots Z_{s-1,s-1}^{}Z_{s,s+1}^{}Z_{s+1,s+2}^{}\dots
	       Z_{r-1,r}^{}Z_{r+1,r+1}^{}\dots Z_{n+1,n+1}^{}&\text{if } r>s,\\
	       Z_{1,1}^{}\dots Z_{r-1,r-1}^{}Z_{r+1,r}^{}Z_{r+2,r+1}^{}\dots
	       Z_{s,s-1}^{}Z_{s+1,s+1}^{}\dots Z_{n+1,n+1}^{} &\text{if } r<s,\\
	       Z_{1,1}^{}\dots Z_{s-1,s-1}^{}Z_{s+1,s+1}^{}Z_{s+2,s+2}^{}\dots
	       Z_{n+1,n+1}^{}  &\text{if } r=s.
		\end{cases}
\end{align}
\end{proposition}
\prf
We will denote the operator $\pi(u_{i,j}(q))$ by
$Y_{i,j}^{\pi}\equiv Y_{i,j}^{\pi}(q)$ in this proof.
From equation~(\ref{eq:comm-ystar}), we get
\begin{align}
  (Y_{r,s}^{\pi})^{*} &=
  (-q)^{s-r}\sum_{\sigma\in\mathscr{S}_{n}}(-q)^{\ell(\sigma)}
  Y_{i_{1},j_{\sigma(1)}}^{\pi}\ldots
  Y_{i_{n},j_{\sigma(n)}}^{\pi}\\
  &= (-q)^{s-r}\sum_{\sigma\in\mathscr{S}_{n}}
  (-q)^{\ell(\sigma)}Y_{1,\widetilde{\sigma}(1)}^{\pi}\ldots
  Y_{r-1,\widetilde{\sigma}(r-1)}^{\pi}
  Y_{r+1,\widetilde{\sigma}(r+1)}^{\pi}\ldots
  Y_{n+1,\widetilde{\sigma}(n+1)}^{\pi},
\end{align}
where $i_{k}$ and $j_{k}$ are as in Lemma~\ref{lmma:perm} and
$\widetilde{\sigma}\in\mathscr{S}_{n+1}$ is the permutation that
takes $i_{k}\mapsto j_{\sigma(k)}$  for all $k\in\{1,2,\ldots,n\}$
and takes $r$ to $s$.

Now let us look at the case $r < s$ first. In this case, one has
\begin{IEEEeqnarray*}{rCl}
	(-q)^{r-s}(Y_{r,s}^{\pi})^{*} &=&
  \sum_{\sigma\in\mathscr{S}_{n}}
  (-q)^{\ell(\sigma)}Y_{1,\widetilde{\sigma}(1)}^{\pi}
  \ldots Y_{r-1,\widetilde{\sigma}(r-1)}^{\pi}
   Y_{r+1,\widetilde{\sigma}(r+1)}^{\pi}
  \ldots  Y_{n+1,\widetilde{\sigma}(n+1)}^{\pi}\\
	&=& Y_{1,1}^{\pi}Y_{2,2}^{\pi}\ldots Y_{r-1,r-1}^{\pi}Y_{r+1,r}^{\pi}
  \ldots Y_{s,s-1}^{\pi}Y_{s+1,s+1}^{\pi}\ldots Y_{n+1,n+1}^{\pi}\\
	& & {}+\sum_{\shortstack{$\scriptstyle \sigma\in\mathscr{S}_{n}$\\
  $\scriptstyle \sigma\neq id$}}
  (-q)^{\ell(\sigma)}Y_{1,\widetilde{\sigma}(1)}^{\pi}\ldots
  Y_{r-1,\widetilde{\sigma}(r-1)}^{\pi}
  Y_{r+1,\widetilde{\sigma}(r+1)}^{\pi}\ldots
  Y_{n+1,\widetilde{\sigma}(n+1)}^{\pi}.
\end{IEEEeqnarray*}
Taking limit as $q\to 0+$, one gets
\begin{align*}
	Z_{r,s}^{*} &=Z_{1,1}^{}\dots Z_{r-1,r-1}^{}Z_{r+1,r}^{}Z_{r+2,r+1}^{}\dots
  Z_{s,s-1}^{}Z_{s+1,s+1}^{}\dots Z_{n+1,n+1}^{}.
\end{align*}

For $r=s$, one has
\begin{IEEEeqnarray*}{rCl}
	(Y_{s,s}^{\pi})^{*} &=& \sum_{\sigma\in\mathscr{S}_{n}}(-q)^{\ell(\sigma)}
  Y_{1,\widetilde{\sigma}(1)}^{\pi}\ldots Y_{s-1,\widetilde{\sigma}(s-1)}^{\pi}
  Y_{s+1,\widetilde{\sigma}(s+1)}^{\pi}\ldots
   Y_{n+1,\widetilde{\sigma}(n+1)}^{\pi}\\
	&=& Y_{1,1}^{\pi}Y_{2,2}^{\pi}\ldots Y_{s-1,s-1}^{\pi}Y_{s+1,s+1}^{\pi}
  \ldots Y_{n+1,n+1}^{\pi}\\
	& & {}+\sum_{\shortstack{$\scriptstyle \sigma\in\mathscr{S}_{n}$\\
                 $\scriptstyle \sigma\neq id$}}
  (-q)^{\ell(\sigma)}Y_{1,\widetilde{\sigma}(1)}^{\pi}\ldots
  Y_{s-1,\widetilde{\sigma}(s-1)}^{\pi}
  Y_{s+1,\widetilde{\sigma}(s+1)}^{\pi}\ldots
   Y_{n+1,\widetilde{\sigma}(n+1)}^{\pi}.
\end{IEEEeqnarray*}
Taking limit as $q\to 0+$, one gets
\begin{align*}
	Z_{s,s}^{*} &=Z_{1,1}^{}\dots Z_{s-1,s-1}^{}Z_{s+1,s+1}^{}Z_{s+2,s+2}^{}
  \dots Z_{n+1,n+1}^{}.
\end{align*}

For $r > s$,  one has
\begin{IEEEeqnarray*}{rCl}
	(Y_{r,s}^{\pi})^{*} &=& (-q)^{s-r}\sum_{\sigma\in\mathscr{S}_{n}}
  (-q)^{\ell(\sigma)}
  Y_{1,\widetilde{\sigma}(1)}^{\pi}\ldots Y_{r-1,\widetilde{\sigma}(r-1)}^{\pi}
  Y_{r+1,\widetilde{\sigma}(r+1)}^{\pi}\ldots
   Y_{n+1,\widetilde{\sigma}(n+1)}^{\pi}\\
	&=& (-q)^{s-r}Y_{1,1}^{\pi}Y_{2,2}^{\pi}\ldots
  Y_{s-1,s-1}^{\pi}Y_{s,s+1}^{\pi}
  \ldots Y_{r-1,r}^{\pi}Y_{r+1,r+1}^{\pi}\ldots Y_{n+1,n+1}^{\pi}\\
	& & {}+(-q)^{s-r}\sum_{\shortstack{$\scriptstyle \sigma\in\mathscr{S}_{n}$\\
   $\scriptstyle \sigma\neq id$}}
  (-q)^{\ell(\sigma)}Y_{1,\widetilde{\sigma}(1)}^{\pi}\ldots
  Y_{r-1,\widetilde{\sigma}(r-1)}^{\pi}
  Y_{r+1,\widetilde{\sigma}(r+1)}^{\pi}\ldots
   Y_{n+1,\widetilde{\sigma}(n+1)}^{\pi}.
\end{IEEEeqnarray*}
Taking limit as $q\to 0+$ and by Lemma~\ref{lmma:perm}, it follows that
\begin{align*}
	Z_{r,s}^{*} &=Z_{1,1}^{}\dots Z_{s-1,s-1}^{}Z_{s,s+1}^{}Z_{s+1,s+2}^{}\dots
  Z_{r-1,r}^{}Z_{r+1,r+1}^{}\dots Z_{n+1,n+1}^{}.
\end{align*}
This completes the proof.
\qed

\begin{proposition}
For any irreducible representation $\pi$ of $A_{n}(q)$, the operators
$Z_{i,j}^{}\equiv Z_{i,j}^{\pi}$ satisfy the following relations:
\begin{equation}
  Z_{i,j}^{*}Z_{r,s}^{}  =  Z_{r,s}^{}Z_{i,j}^{*}, \qquad i\neq r,\; j\neq s.
  \label{eq:comm-y9}
\end{equation}
\end{proposition}
\prf
This is immediate from (\ref{eq:comm-t7}) and Proposition~\ref{prop:key-prop}.
\qed

  \subsection{The crystallized $C^{*}$-algebra $A_{n}(0)$}
\begin{theorem}
There is a universal $C^{*}$-algebra generated by elements $z_{i,j}^{}$, $1\leq
i,j\leq n+1$ satisfying the following relations:
\begin{IEEEeqnarray}{rClll}
   z_{i,j}^{}z_{i,l}^{} &=& 0 & \text{if }
            &j<l,\label{eq:comm-z1}\\
   z_{i,j}^{}z_{k,j}^{} &=& 0 & \text{if }
            &i<k,\label{eq:comm-z2}\\
            &&&&\nonumber\\
    z_{i,l}^{}z_{k,j}^{} - z_{k,j}^{}z_{i,l}^{} &=& 0 \quad
     & \text{if } &i<k\text{ and }j<l.\label{eq:comm-z4}\\
     &&&&\nonumber\\
     z_{i,l}^{}z_{k,j}^{}  &=&0 & \text{if }& i<k,\,j<l \text{ and }
           \max\{i,j\}\geq \min\{k,l\},\label{eq:comm-z3}\\
      &&&&\nonumber\\
  z_{i,j}^{}z_{k,l}^{} - z_{k,l}^{}z_{i,j}^{}  &=& 0 &\text{if }&i<k,\,j<l \text{ and }
           \max\{i,j\}+1 < \min\{k,l\},\label{eq:comm-z5}\\
    &&&&\nonumber\\
z_{i,j}^{}z_{k,l}^{} - z_{k,l}^{}z_{i,j}^{}&=& z_{i,l}^{}z_{k,j}^{}\quad
     & \text{if }&i<k,\,j<l \text{ and }
           \max\{i,j\}+1= \min\{k,l\},\label{eq:comm-z6}\\
    &&&&\nonumber\\
    z_{1,1}^{}z_{2,2}^{}\ldots z_{n+1,n+1}^{}&=& 1,\label{eq:comm-z7a}\\
     &&&&\nonumber\\
     \IEEEeqnarraymulticol{5}{l}{\hspace{-3em}z_{r,s}^{*}  =
        \begin{cases}
         \left(z_{1,1}^{}\dots z_{s-1,s-1}^{}\right)\left(z_{s,s+1}^{}z_{s+1,s+2}^{}\dots
         z_{r-1,r}^{}\right)\left(z_{r+1,r+1}^{}\dots z_{n+1,n+1}^{}\right)&\text{if } r>s,\\
       \left(z_{1,1}^{}\dots z_{r-1,r-1}^{}\right)\left(z_{r+1,r}^{}z_{r+2,r+1}^{}\dots
       z_{s,s-1}^{}\right)\left(z_{s+1,s+1}^{}\dots z_{n+1,n+1}^{}\right) &\text{if } r<s,\\
      \left(z_{1,1}^{}\dots z_{s-1,s-1}^{}\right)\left(z_{s+1,s+1}^{}z_{s+2,s+2}^{}\dots
       z_{n+1,n+1}^{}\right)  &\text{if } r=s.
    \end{cases}}\nonumber\\
    &&&&\label{eq:comm-z7}\\
  \IEEEeqnarraymulticol{5}{l}{z_{i,j}^{*}z_{r,s}^{}  =  z_{r,s}^{}z_{i,j}^{*}, \qquad i\neq r,\; j\neq s.} \label{eq:comm-z8} 
\end{IEEEeqnarray}
\end{theorem}
\prf
Let $\mathcal{H}$ be the Hilbert space $\ell^{2}(\mathbb{N})$.
Let $1\leq k\leq n$ and let
\[
\pi(z_{i,j}^{})=
    \begin{cases}
        S &\text{if }i=j=k,\cr
        S^{*} & \text{if } i=j=k+1,\cr
        |e_{0}\rangle\,\langle e_{0}| &  \text{if } i=k,\;\; j=k+1,\cr
        |e_{0}\rangle\,\langle e_{0}|   &
        \text{if }i=k+1,\;\; j=k,\cr
        \delta_{i,j}I & \text{otherwise},
		\end{cases}
\]
where $|e_{0}\rangle\,\langle e_{0}|$ is the projection onto the span of
$e_{0}$. Then $\pi$ gives a representation of the relations
(\ref{eq:comm-z1}--\ref{eq:comm-z8}). 

It follows from the relations (\ref{eq:comm-z5}--\ref{eq:comm-z8}) that
\begin{IEEEeqnarray}{rClCl}
 \sum_{k=1}^{j}z_{j,k}^{}z_{j,k}^{*} & = & 1 &=&\sum_{k=j}^{n+1}z_{k,j}^{*}z_{k,j}^{},\qquad 1\leq i\leq n+1.\label{eq:comm-z9}
\end{IEEEeqnarray}
Therefore we have $\|\pi(z_{i,j}^{})\|\leq 1$ for all $i\geq j$.
Using (\ref{eq:comm-z7}), one then concludes that 
$\|\pi(z_{i,j}^{})\|\leq 1$ for all $i<j$. Thus
$\|\pi(z_{i,j}^{})\|\leq 1$ for all $i,j$. 
The rest of the proof is now standard and follows from the fact that for any polynomial $a$ in the $z_{i,j}^{}$'s and $z_{i,j}^{*}$'s, there is a positive real $K$ such that $\|\pi(a)\|\leq K$ for any representation $\pi$ of the relations
(\ref{eq:comm-z1}--\ref{eq:comm-z8}).
\qed

\begin{remark}
The above definition of crystallization relies on the observations that \\[1ex]
1.
  the limits $\lim_{q\to 0+}(-q)^{\min\{i-j,0\}}\pi_{i,j}(u_{i,j}(q))$ exist for each irreducible representation $\pi$, \\
  2.
  they lead to the same set of relations, and\\
  3.
  these relations ensure that for each representation of them on a Hilbert space, the norm of each generating element is bounded by a fixed constant.\\[1ex]
The knowledge of all irreducible representations of $C(G_{q})$, due to
Soibelman, plays a very crucial role in the first two points above, while
unitarity of the matrix $(\!(u_{i,j})\!)$ ensures point 3 above. Soibelman's
results remain valid for deformations of compact Lie groups for all the other
types ($B$, $C$, $D$, $E$, $F$ and $G$). Using his result, it
is possible to carry out steps 1 and 2 above for other types as well
and unitarity of $(\!(u_{i,j})\!)$ will then guarantee that step 3 remains
valid. Thus one can define crystallization of $C(G_{q})$ in
the remaining cases as well following the above mechanism.
However, as has
been mentioned in the  introduction, our aim in this paper is to study
the irreducible representations of the crystallized $C^{*}$-algebra
in the type $A_{n}$ case and accordingly we will focus on irreducible representations
of the crystallized algebra $A_{2}(0)$ in the rest of this paper, and plan to take 
up the irreducible representations of $A_{n}(0)$ for $n>2$ in a separate article.
\end{remark}
\begin{remark}
After the first version of this paper came out on the arXiv, Matassa \& Yuncken
(\cite{MatYun-2022aa}) introduced crystallizations of the $C^*$-algebras $C(G_
{q})$ for any compact semisimple Lie group $G$, i.e.\ for a much broader class of
$q$-deformations than considered here. They also go on to prove that the crystallized
$C^{*}$-algebras are higher rank graph $C^{*}$-algebras. However, the trade-off is that their
construction is more involved and uses results on crystal bases due to
Kashiwara et al. Our approach here is simpler and we get the crystallized algebra as
a $C^{*}$-algebra given by a finite set of generators and relations, which
make it convenient to study certain properties, for example, finding its irreducible representations
that we take up in the rest of the paper in the rank two case.
\end{remark}

 \section{The case $n=2$}
In this section we will focus our attention to the case $n=2$ and study 
the irreducible representations of the $C^*$-algebra $A_{2}(0)$.
\subsection{Irreducible representations}
Let us start by introducing a family of representations of $A_n(0)$ that 
are analogous
to the irreducible representations of $A_n(q)$ for $q\neq 0$. As mentioned in
Section~2, Theorem~6.2.7 in \cite{KorSoi-1998ab} tells us that the family
$\psi_{\lambda,\omega}^{(q)}$ where $\lambda\in (S^{1})^{n}$ and
$\omega$ is a reduced word in $\mathscr{S}_{n+1}$, gives all the ireducible
representations of the $C^{*}$-algebra $A_{n}(q)$ for $q\in (0,1)$. Let us
define
\begin{IEEEeqnarray}{rCl}
	Z_{i,j}^{} &=& \lim_{q\to 0+}(-q)^{\min\{i-j,0\}}
	 \psi_{\lambda,\omega}^{(q)}(u_{i,j}^{}(q)),
	   \label{eq:irred-q0-1}
\end{IEEEeqnarray}
where $\psi_{\lambda,\omega}^{(q)}$ is the representation of $A_{n}(q)$ defined
in Section~2. Note that by Proposition~\ref{prop:key-prop}, the above limits
exist and the operators $Z_{i,j}^{}$ obey the relations
(\ref{eq:comm-z1}--\ref{eq:comm-z8}). Therefore
\begin{IEEEeqnarray}{rCl}
\psi_{\lambda,\omega}^{}(z_{i,j}^{})=
 \lim_{q\to 0+}(-q)^{\min\{i-j,0\}}
 \psi_{\lambda,\omega}^{(q)}(u_{i,j}^{}(q)),\qquad
  i,j\in\{1,2,\ldots,n+1\}\label{eq:irred-q0-2}
\end{IEEEeqnarray}
defines a representation of $A_{n}(0)$ on the Hilbert space
$\ell^{2}(\mathbb{N}^{\ell(\omega)})$, where $\ell(\omega)$ denotes the length
of the element $\omega$.

\begin{theorem}\label{thm:irred-q0-3}
For any $(\lambda,\mu)\in (S^{1})^{2}$ and a reduced word $\omega$ in
$\mathscr{S}_{3}$, the representation $ \psi_{(\lambda,\mu),\omega}$ given by
(\ref{eq:irred-q0-2}) is irreducible.

If $((\lambda,\mu),\omega)\neq ((\lambda',\mu'),\omega')$, then the
representations $\psi_{(\lambda,\mu),\omega}$ and
$\psi_{(\lambda',\mu'),\omega'}$ are inequivalent.
\end{theorem}
\prf
We will denote the rank one projection $|e_{0}\rangle\,\langle e_{0}|$ by
$P_{0}$. For $\omega=id$, the representation $ \psi_{(\lambda,\mu),\omega}$ is
one dimensional and hence irreducible. Let $\omega=s_{1}$. In this case, $
\psi_{(\lambda,\mu),\omega}$ is given by
\begin{IEEEeqnarray*}{rClrClrCl}
  z_{1,1}^{}&\mapsto & \lambda S, \qquad &
    z_{1,2}^{}&\mapsto & \lambda P_{0}, \qquad &
      z_{1,3}^{}&\mapsto & 0, \qquad \\
        z_{2,1}^{}&\mapsto & \bar{\lambda}\mu P_{0}, \qquad &
          z_{2,2}^{}&\mapsto & \bar{\lambda}\mu S^{*}, \qquad &
            z_{2,3}^{}&\mapsto & 0, \qquad \\
              z_{3,1}^{}&\mapsto & 0, \qquad &
                z_{3,2}^{}&\mapsto & 0, \qquad &
                  z_{3,3}^{}&\mapsto & \bar{\mu}I .
\end{IEEEeqnarray*}
Therefore for any $j,k\in\mathbb{N}$, one has
\[
|e_{j}\rangle\,\langle e_{k}|=
  \lambda^{j-k-1}(z_{1,1}^{*})^{j}z_{1,2}^{}z_{1,1}^{k}.
\]
Therefore $\mathcal{K}(\ell^{2}(\mathbb{N}))\subseteq
\psi_{(\lambda,\mu),\omega}(A_{2}(0))$. This imples
$\psi_{(\lambda,\mu),\omega}$ is irreducible. Proof for $\omega=s_{2}$ is
similar.

Next, let us take $\omega=s_{1}s_{2}$. In this case we have
\begin{IEEEeqnarray*}{rClrClrCl}
  z_{1,1}^{}&\mapsto & \lambda S\otimes I, \qquad &
    z_{1,2}^{}&\mapsto & \lambda P_{0}\otimes S, \qquad &
      z_{1,3}^{}&\mapsto & \lambda P_{0}\otimes P_{0}, \qquad \\
        z_{2,1}^{}&\mapsto & \bar{\lambda}\mu P_{0}\otimes I, \qquad &
          z_{2,2}^{}&\mapsto & \bar{\lambda}\mu S^{*}\otimes S, \qquad &
  z_{2,3}^{}&\mapsto & \bar{\lambda}\mu S^{*}\otimes P_{0}, \qquad \\
              z_{3,1}^{}&\mapsto & 0, \qquad &
                z_{3,2}^{}&\mapsto & \bar{\mu}I\otimes P_{0}, \qquad &
                  z_{3,3}^{}&\mapsto & \bar{\mu}I\otimes S^{*} .
\end{IEEEeqnarray*}
In this case, for any $j,k\in\mathbb{N}$, we have
\begin{IEEEeqnarray*}{rCl}
|e_{j}\rangle\,\langle e_{k}|\otimes I &=&
  \lambda^{j-k-1}(z_{1,1}^{*})^{j}z_{1,2}^{}z_{1,1}^{k},\\
  I\otimes |e_{j}\rangle\,\langle e_{k}| &=&
    \mu^{j-k+1}z_{3,3}^{j}z_{3,2}^{}(z_{3,3}^{*})^{k}.
\end{IEEEeqnarray*}
Thus
$\mathcal{K}(\ell^{2}(\mathbb{N})\otimes \ell^{2}(\mathbb{N}))\subseteq
\psi_{(\lambda,\mu),\omega}(A_{2}(0))$, which imples
$\psi_{(\lambda,\mu),\omega}$ is irreducible. Proof for
$\omega=s_{2}s_{1}$ is similar.

Finally take $\omega=s_{1}s_{2}s_{1}$. Here one has
\begin{IEEEeqnarray*}{rClrClrCl}
  z_{1,1}^{}&\mapsto & \lambda S\otimes I\otimes S, \qquad &
    z_{1,2}^{}&\mapsto & \lambda P_{0}\otimes S\otimes S^{*}
      +S\otimes I\otimes P_{0}, \qquad &
      z_{1,3}^{}&\mapsto & \lambda P_{0}\otimes P_{0}\otimes I, \qquad \\
        z_{2,1}^{}&\mapsto & \bar{\lambda}\mu
          P_{0}\otimes I\otimes S, \qquad &
          z_{2,2}^{}&\mapsto & \bar{\lambda}\mu
          S^{*}\otimes S\otimes S^{*}, \qquad &
  z_{2,3}^{}&\mapsto & \bar{\lambda}\mu S^{*}\otimes P_{0}\otimes I, \qquad \\
    z_{3,1}^{}&\mapsto & I\otimes P_{0}\otimes P_{0}, \qquad &
    z_{3,2}^{}&\mapsto & \bar{\mu}I\otimes P_{0}\otimes S^{*}, \qquad &
    z_{3,3}^{}&\mapsto & \bar{\mu}I\otimes S^{*}\otimes I .
\end{IEEEeqnarray*}
Now observe that for any $j,k,m,n\in\mathbb{N}$,
\begin{IEEEeqnarray*}{rCl}
|e_{j}\rangle\,\langle e_{k}|\otimes |e_{m}\rangle\,\langle e_{n}|\otimes I &=&
  c_{1}z_{3,3}^{m}
  \left(z_{2,3}^{j}z_{1,3}^{}(z_{2,3}^{*})^{k} \right)(z_{3,3}^{*})^{n},\\
  I\otimes |e_{j}\rangle\,\langle e_{k}|\otimes |e_{m}\rangle\,\langle e_{n}|
  &=&
  c_{2}z_{3,3}^{j}
  \left(z_{3,2}^{m}z_{3,1}^{}(z_{3,2}^{*})^{n} \right)(z_{3,3}^{*})^{k},
\end{IEEEeqnarray*}
where $c_{1}$ and $c_{2}$ are two scalars involving $\lambda$ and $\mu$.
Therefore it follows that
\[
\mathcal{K}(\ell^{2}(\mathbb{N})\otimes
\ell^{2}(\mathbb{N})\otimes \ell^{2}(\mathbb{N}))\subseteq
\psi_{(\lambda,\mu),\omega}(A_{2}(0)),
\]
which imples that $\psi_{(\lambda,\mu),\omega}$ is irreducible.

The inequivalence of different $\psi_{(\lambda,\mu),\omega}$'s follow by
comparing the spectrum of the operators
$\psi_{(\lambda,\mu),\omega}(z_{i,j}^{})$'s.
\qed

The next theorem states that the ones listed in the previous theorem are in fact all the irreducible representations of $A_{2}(0)$. The proof is computation-heavy and has been given in the last section. 
\begin{theorem}\label{thm:irred-q0-4}
Let $\pi$ be an irreducible representation of $A_{2}(0)$ on a Hilbert space $\mathcal{H}$. Then there exist $\lambda,\mu\in S^{1}$ and a reduced word 
$\omega\in \mathscr{S}_{3}$ such that $\pi$ is equivalent to 
$\psi_{(\lambda,\mu),\omega}$.
\end{theorem}

\begin{remark}
Theorem~\ref{thm:irred-q0-3} and Theorem~\ref{thm:irred-q0-4} together gives us
the family of all inequivalent irreducible representations of the
$C^{*}$-algebra $A_{2}(0)$. They are naturally parametrized exactly as their
counterparts for nonzero $q$, and in fact they arise as limits of
the irreducible representations of the $C^{*}$-algebras $A_{2}(q)$ as $q\to 0+$.
\end{remark}

\begin{remark}
In the course of the above proof, we have also observed that for any irreducible
representation $\pi$ of $A_{2}(0)$ on a Hilbert space $\mathcal{H}_{\pi}$, one
has $\mathcal{K}(\mathcal{H}_{\pi})\subseteq \pi(A_{2}(0))$. Thus $A_{2}(0)$ is
a type I $C^{*}$-algebra.
\end{remark}

\subsection{Applications}
In this subsection, we derive a few properties of the $C^{*}$-algebra $A_{2}(0)$
using the results in the previous subsection.

For the $C^{*}$-algebra $A_{1}(0)$, each infinite dimensional irreducible
representation $\pi$ acts on the Hilbert space $\ell^{2}(\mathbb{N})$ and one
has $\pi(A_{1}(0))=\mathscr{T}$, where $\mathscr{T}$ denotes the Toeplitz
algebra. In the case of $A_{2}(0)$, each infinite dimensional irreducible $\pi$
acts on a Hilbert space of the form $\ell^{2}(\mathbb{N})^{\otimes k}$, and it
is straightforward to verify that $\pi(A_{2}(0))\subseteq
\mathscr{T}(\ell^{2}(\mathbb{N})^{\otimes k}):=\mathscr{T}^{\otimes k}$. It
would be natural to ask whether
$\pi(A_{2}(0))=\mathscr{T}(\ell^{2}(\mathbb{N})^{\otimes k})$. The next
proposition says that it is not the case.

\begin{proposition}\label{prop:image-toeplitz}
Let $\omega=s_{1}s_{[1,2]}$. Then
$\psi_{\omega}(A_{2}(0))\subsetneqq \mathscr{T}^{\otimes 3}$.
\end{proposition}
\prf
We have already observed that $\psi_{\omega}(A_{2}(0))\subseteq
\mathscr{T}^{\otimes 3}$. We will produce an element in $\mathscr{T}^{\otimes
3}$ that does not belong to the image $\psi_{\omega}(A_{2}(0))$. For $\lambda
\in S^{1}$, let $\sigma_{\lambda}$ denote the $C^*$-homomorphism from
$\mathscr{T}$ to $\mathbb{C}$ obtained by composing the evaluation map
$ev_{\lambda }:C(S^{1})\to\mathbb{C}$ with the canonical projection map from
$\mathscr{T}$ to $C(S^{1})\cong \mathscr{T}/\mathcal{K}$. Let $\phi_{\lambda
}:\mathscr{T}^{\otimes 3}\to\mathbb{C}$ be the $C^*$-homomorphism given by
$\phi_{\lambda }=\sigma_{\lambda }\otimes \sigma_{1}\otimes \sigma_{\bar{\lambda
}}$. Then
\[
	\phi_{\lambda }\left(\psi_{(1,1),\omega}(z_{i,j}^{})\right) =
  \sum_{k, l}\sigma_{\lambda}\left(\psi_{(1,1),s_{1}}
  (z_{i,k}^{})\right)\sigma_{1}\left(\psi_{(1,1),s_{2}}
  (z_{k,l}^{})\right)
	    \sigma_{\bar{\lambda }}\left(\psi_{(1,1),s_{1}}
      (z_{l,j}^{})\right) =\delta_{i,j}.
\]
Take $S\otimes S\otimes S^{*}\in\mathscr{T}^{\otimes 3}$. Assume $S\otimes
S\otimes S^{*}\in \psi_{(1,1),\omega}(A_{2}(0))$. Take $\epsilon<\frac{1}{2}$.
Then there exists a polynomial $p\equiv p(z_{i,j}^{}, z_{i,j}^{*}: i,j)$ in the
$z_{i,j}^{}$'s and $z_{i,j}^{*}$'s such that
\begin{align*}
	\|S\otimes S\otimes S^{*}-\psi_{(1,1),\omega}(p)\| &< \epsilon.
\end{align*}
Now $\phi_{\lambda }(S\otimes S\otimes S^{*})=\lambda ^{2}$ and $\phi_{\lambda
}(\psi_{(1,1),\omega}(p))=\sum c_{\gamma}$, where $c_{\gamma}$ are the
coefficients of the monomials in $z_{i,i}^{}$ and $z_{j,j}^{*}$'s. Therefore for
any $\lambda $ and $\lambda '$ in $S^{1}$, we have
\begin{IEEEeqnarray*}{rClrCl}
	|\lambda ^{2}-\sum c_{\gamma}| &<& \epsilon,\qquad &
	|(\lambda ')^{2}-\sum c_{\gamma}| &< &\epsilon.
\end{IEEEeqnarray*}
This implies $|\lambda ^{2}-(\lambda')^{2}|<2\epsilon$ for all
$\lambda ,\lambda'\in S^{1}$, which is false.
\qed

Let $\omega$ be as above. Since $\psi_{(\lambda,\mu),\omega}(A_{2}(0))\subseteq
\mathscr{T}^{\otimes 3}$, one can view the irreducible representation
$\psi_{(\lambda,\mu),\omega}$ as a $C^{*}$-homomorphism into the $C^{*}$-algebra
$\mathscr{T}^{\otimes 3}$. Let $u\in C(S^{1})$ denote the function
$\lambda\mapsto \lambda$ and $\one\in C(S^{1})$ denote the function
$\lambda\mapsto 1$. Let us define a map $\phi$ from $A_{2}(0)$ to the
$C^{*}$-algebra $C(S^{1})\otimes C(S^{1})\otimes \mathscr{T}^{\otimes 3}$ as
follows:
\begin{IEEEeqnarray*}{rCl}
  z_{i,j}^{} &\mapsto &
  \begin{cases}
    u\otimes \one \otimes \psi_{(1,1),\omega}(z_{i,j}^{}) & \text{if }i=1,\\
    u^{*}\otimes u \otimes \psi_{(1,1),\omega}(z_{i,j}^{}) & \text{if }i=2,\\
    \one \otimes u^{*} \otimes \psi_{(1,1),\omega}(z_{i,j}^{}) & \text{if }i=3.
  \end{cases}
\end{IEEEeqnarray*}
\begin{proposition}
The map $\phi$ defined above extends to an injective $*$-homomorphism from
$A_{2}(0)$ to $C(S^{1})\otimes C(S^{1})\otimes \mathscr{T}^{\otimes 3}$.
\end{proposition}
\prf
Since the elements $\phi(z_{i,j}^{})$ obey the relations
(\ref{eq:comm-z1}--\ref{eq:comm-z8}), the map $\phi$ extends to a
$C^{*}$-homomorphism from $A_{2}(0)$ into $C(S^{1})\otimes C(S^{1})\otimes
\mathscr{T}^{\otimes 3}$. For injectivity, note that
$\psi_{(\lambda,\mu),\tau}$, where $\lambda,\mu\in S^{1}$ and $\tau$ is a
reduced word in $\mathscr{S}_{3}$, constitute all irreducible representations of
$A_{2}(0)$. Therefore it is enough to show that each $\psi_{(\lambda,\mu),\tau}$
factorizes through $\phi$. Let $\sigma_{\lambda}:\mathscr{T}\to\mathbb{C}$ be as
in the proof of Proposition~\ref{prop:image-toeplitz}. Define
\begin{IEEEeqnarray*}{rCl}
  \theta_{\tau} &=&
   \begin{cases}
   \sigma_{1}\otimes \sigma_{1} \otimes \sigma_{1} & \text{if }\tau=id,\\
     \id\otimes \sigma_{1}\otimes \sigma_{1} & \text{if }\tau=s_{1},\\
     \sigma_{1}\otimes \id \otimes \sigma_{1} & \text{if }\tau=s_{2},\\
     \id\otimes \id \otimes \sigma_{1} & \text{if }\tau=s_{1}s_{2},\\
     \sigma_{1}\otimes \id \otimes \id & \text{if }\tau=s_{2}s_{1},\\
     \id\otimes \id \otimes \id & \text{if }\tau=s_{1}s_{2}s_{1}.
   \end{cases}
\end{IEEEeqnarray*}
Then
$\psi_{(\lambda,\mu),\tau}=(ev_{\lambda}\otimes
ev_{\mu}\otimes\theta_{\tau})\circ\phi$.
\qed

Through the embeddings $C(S^{1})\hookrightarrow
\mathcal{L}(\ell^{2}(\mathbb{Z}))$ and $\mathscr{T}^{\otimes 3}\hookrightarrow
\mathcal{L}(\ell^{2}(\mathbb{N})^{\otimes 3})$, the above map $\phi$ gives a
faithful representation of $A_{2}(0)$ acting on the Hilbert space
$\ell^{2}(\mathbb{Z})\otimes \ell^{2}(\mathbb{Z})\otimes
\ell^{2}(\mathbb{N})\otimes \ell^{2}(\mathbb{N})\otimes \ell^{2}(\mathbb{N})$.
This is the analogue of the representation $\pi_{0}$ of $A_{1}(0)$ used in
\cite{ChaPal-2022ey} by Chakraborty \& Pal.

We prove next that  $A_{2}(0)$ admits a coproduct that makes it a
compact quantum semigroup. However, it does not make it a compact quantum group.
\begin{proposition}
The map $\Delta:A_{2}(0)\to A_{2}(0) \otimes A_{2}(0)$ given on the generators
by
\[
\Delta(z_{i,j}^{})=\sum_{k=\min\{i,j\}}^{\max\{i,j\}}z_{i,k}^{}\otimes z_{k,j}^{}
\]
together with $\epsilon=\psi_{(1,1),id}$ from $A_{2}(0)$ to $\mathbb{C}$ makes
$(A_{2}(0),\Delta,\epsilon)$ a $C^{*}$-bialgebra.
\end{proposition}
\prf
Note that since $A_{2}(0)$ is of type I, the tensor product $C^{*}$-algebra
$A_{2}(0)\otimes A_{2}(0)$ is unique.
In order that $\Delta$ defined on the generating elements as above extend to a
$C^{*}$-homomorphism from $A_{2}(0)$ to  $A_{2}(0)\otimes A_{2}(0)$, it is
enough to verify that the elements $\Delta(z_{i,j}^{})$ obey the same relations
(\ref{eq:comm-z1}--\ref{eq:comm-z8}) as the $z_{i,j}^{}$'s. Recall from
Section~2 that for any irreducible representation $\pi$ of $A_{n}(q)$, the
operators $Y_{i,j}^{}=\lim_{q\to 0+}(-q)^{\min\{i-j,0\}}\pi(u_{i,j}(q))$ satisfy
these same relations.

Let $q\in(0,1)$. Since $\Delta_{q}$ is a $C^{*}$-homomorphism from $A_{2}(q)$ to
$A_{2}(q)\otimes A_{2}(q)$, the elements
$\Delta_{q}(u_{i,j}^{})=\sum_{k=1}^{3}u_{i,k}^{}\otimes u_{k,j}^{}$ obey the
same commutation relations as the $u_{i,j}^{}$'s. Hence for $\lambda,\mu\in
(S^{1})^{2}$ and $\omega,\sigma$ reduced words in $\mathscr{S}_{3}$, the
elements
\[
\left(\psi_{\lambda,\omega}^{(q)}\otimes \psi_{\mu,\sigma}^{(q)}\right)
 \left(\Delta_{q}(u_{i,j}^{})\right)
\]
also satisfy the same relations.
Now observe that
\[
\lim_{q\to 0+}
\left(\psi_{\lambda,\omega}^{(q)}\otimes \psi_{\mu,\sigma}^{(q)}\right)
\left(\Delta_{q}((-q)^{\min\{i-j,0\}}u_{i,j}^{})\right)=
\left(\psi_{\lambda,\omega}\otimes \psi_{\mu,\sigma}\right)
 \left(\sum_{k=\min\{i,j\}}^{\max\{i,j\}}z_{i,k}^{}\otimes z_{k,j}^{}\right).
\]
Since $A_{2}(0)$ is of type I, the representations
$\psi_{\lambda,\omega}\otimes \psi_{\mu,\sigma}$
give all irreducible representations of the $C^{*}$-algebra
$A_{2}(0)\otimes A_{2}(0)$. Therefore it follows that the $\Delta(z_{i,j}^{})$'s
obey the relations (\ref{eq:comm-z1}--\ref{eq:comm-z8}).

The identities $(\Delta\otimes \id)\Delta=(\id\otimes\Delta)\Delta$ and
$(\epsilon\otimes \id)\Delta=(\id\otimes\epsilon)\Delta=\id$ follow by
evaluating both sides on the generating elements.
\qed

\begin{remark}
The proof of existence of the homomorphism $\Delta$ above requires proving that
the elements $\sum_{k=\min\{i,j\}}^{\max\{i,j\}}z_{i,k}^{}\otimes z_{k,j}^{}$ in the
tensor product $A_{2}(0)\otimes A_{2}(0)$ satisfy the defining relations of the
generators of $A_{2}(0)$. While this should be a direct algebraic verification,
a direct verification has eluded us. Therefore we have used the knowledge
of the irreducible representations to give an indirect proof of that fact here.
However, we believe that  a brute force verification should in principle be
possible for all $n$.
\end{remark}

We will next prove that $(A_{2}(0),\Delta)$ is not a compact quantum group.
Recall (\cite{Wor-1998aa}) that a separable unital $C^{*}$-algebra $A$ together
with a unital $C^{*}$-homomorphism $\Delta:A\to A\otimes A$ is called a compact
quantum group if the following two conditions are satisfied:
\begin{enumerate}
\item
$(\Delta\otimes \id)\Delta=(\id\otimes\Delta)\Delta$.
\item
The sets $\left\{(1\otimes a)\Delta(b): a,b\in A\right\}$ and
$\left\{(a\otimes 1)\Delta(b): a,b\in A\right\}$ are total in $A\otimes A$.
\end{enumerate}
We will show that the second condition above is violated for $A_{2}(0)$.

Let us denote the generators of $A_1(0)$ by $y_{i,j}^{}$. Then one can
directly verify that the elements $\sum_{k=\min\{i,j\}}^{\max\{i,j\}}y_{i,k}^{}\otimes y_{k,j}^{}$ in the tensor product $A_1(0) \otimes
A_1(0)$  satisfy the same relations as the $y_{i,j}^{}$'s and hence
$\Delta^{(1)}: A_1(0) \to A_1(0) \otimes A_1(0)$ given by
$y_{i,j}^{}\mapsto \sum_{k=\min\{i,j\}}^{\max\{i,j\}}y_{i,k}^{}\otimes y_{k,j}^{}$
makes $(A_{1}(0),\Delta^{(1)})$ a compact quantum semigroup. Now let
\[
w_{i,j}^{}=
\begin{cases}
  y_{i,j}^{} &\text{if }1\leq i,j\leq 2,\\
  \delta_{i,j} & \text{if $i=3$ or $j=3$}.
\end{cases}
\]
Then one can verify that the $w_{i,j}^{}$'s obey the relations
(\ref{eq:comm-z1}--\ref{eq:comm-z8}). Therefore the map $\phi:A_{2}(0)\to
A_{1}(0)$ given by
\[
\phi(z_{i,j}^{})=w_{i,j}^{}
\]
extends to a surjective $C^{*}$-homomorphism and it satisfies
$(\phi\otimes\phi)\Delta=\Delta^{(1)}\phi$. Thus $(A_{1}(0),\Delta^{(1)})$ is a
quantum subsemigroup of $(A_{2}(0),\Delta)$.

Next, note  that $A_{1}(0)$ is not a compact quantum group with the above
coproduct. The proof is simple, but since it is hard to find a reference, let us
give a quick proof here. Let us first recall from Section~3 that the
$C^{*}$-algebra $A_{1}(0)$ is given by the generators $y_{i,j}^{}$ satisfying
the relations
\begin{IEEEeqnarray*}{rClrCl}
y_{1,1}^{*}y_{1,1}^{}+y_{2,1}^{*}y_{2,1}^{} &=& 1,\qquad
       & y_{1,1}^{}y_{1,1}^{*} &=& 1,\\
      y_{1,1}^{}y_{1,2}^{} &=& 0,
    &y_{1,1}^{}y_{2,1}^{} &=& 0,\yesnumber\label{eq:su02-relations}\\
y_{2,1}^{}y_{1,2}^{}&=&y_{1,2}^{}y_{2,1}^{},\quad &&&\\
y_{1,1}^{*}&=& y_{2,2}^{},& y_{2,1}^{*}&=& y_{1,2}^{}.
\end{IEEEeqnarray*}
By the results of Woronowicz (\cite{Wor-1998aa}), it follows that if
$(A_{1}(0),\Delta^{(1)})$ is a compact quantum group, then $A_{1}(0)$ has a
dense subalgebra $B$ containing the $y_{i,j}^{}$'s such that
$\Delta^{(1)}(B)\subseteq B\otimes_{alg} B$ and there is an antipode map $S:B\to
B$ such that
\[
m(\id\otimes S)\Delta^{(1)} (b)=\epsilon(b)\cdot 1=
m(S\otimes \id)\Delta^{(1)} (b),\qquad\text{for all }b\in B,
\]
where $m:B\otimes_{alg} B \to B$ is the multiplication map $a\otimes b\mapsto
ab$ and $\epsilon:B\to\mathbb{C}$ is the counit given by $y_{i,j}^{}\mapsto
\delta_{i,j}$. Therefore taking $b=y_{1,1}^{}$, we get
\[
y_{1,1}^{}S(y_{1,1}^{})= 1= S(y_{1,1}^{})y_{1,1}^{}.
\]
From the equality $1=S(y_{1,1}^{})y_{1,1}^{}$ and
$y_{1,1}^{}(y_{1,1}^{})^{*}=1$, we get
$S(y_{1,1}^{})=(y_{1,1}^{})^{*}$. But then
$S(y_{1,1}^{})y_{1,1}^{}=(y_{1,1}^{})^{*}y_{1,1}^{}
=1-(y_{2,1}^{})^{*}y_{2,1}^{}\neq 1$. Thus we get a contradiction. Thus
$(A_{1}(0),\Delta^{(1)})$ is not a compact quantum group.

Finally, assume that the density condition in the definition of a compact
quantum group holds for $(A_{2}(0),\Delta)$. Then applying the homomorphism
$\phi\otimes\phi$ and using the fact that it is surjective, it follows that
$\{(1\otimes a)\Delta^{(1)}(b): a,b\in A_{1}(0)\}$ and $\{(a\otimes
1)\Delta^{(1)}(b): a,b\in A_{1}(0)\}$ are total in $A_{1}(0)\otimes A_{1}(0)$.
But this implies that $(A_{1}(0),\Delta^{(1)})$ is a compact quantum group
which is false. Thus $(A_{2}(0),\Delta)$ is not a compact quantum group.

\section{Proof of Theorem~\ref{thm:irred-q0-4}}
We will give a proof of Theorem~\ref{thm:irred-q0-4} in this section.
We first prove a few lemmas on the generating elements and certain products of them that will be needed in the proof.
\subsection{Partial isometries}
We will prove in this subsection that the generating elements and certain
products of them  are  partial isometries. This will be needed in the
computations that follow in obtaining all the irreducible representations.

Define $p_{i,j}^{}=z_{i,j}^{*}z_{i,j}^{}$ and
$q_{i,j}^{}=z_{i,j}^{}z_{i,j}^{*}$. In terms of the $p_{i,j}^{}$'s and
$q_{i,j}^{}$'s, the relations (\ref{eq:comm-z9}) can be written as
\begin{equation}
\begin{aligned}
p_{1,1}^{}+p_{2,1}^{} + p_{3,1}^{} &= 1,\qquad
   & q_{3,1}^{} +q_{3,2}^{}+q_{3,3}^{} &= 1,\\
p_{2,2}^{}+p_{3,2}^{} &= 1, & q_{2,1}^{} + q_{2,2}^{} &= 1,\\
p_{3,3}^{}&= 1, & q_{1,1}^{} &= 1,
\end{aligned} \label{eq:comm-unitary1a}
\end{equation}
We will first show that the generators $z_{i,j}^{}$ are all partial isometries.
For this, we will use the following result.
\begin{lemma}[Erdelyi \cite{Erd-1968aa}, Halmos-Wallen
    \cite{HalWal-1970aa}]\label{lemma:partial-iso-1}
  Let $ u $ and $ v $ be partial isometries. Then $ uv $ is a partial isometry
    if and only if $ u^{*}u $ and $ vv^{*} $ commute.
\end{lemma}

\begin{lemma}\label{lemma:partial-iso-2}
  The family $ \{p_{i,j}^{}, q_{i,j}^{}: 1\leq j \leq i \leq 3\}$ is a commuting
  family of projections.
\end{lemma}
\prf
First observe that
\[
p_{1,1}^{}=z_{1,1}^{*}z_{1,1}^{}=z_{2,2}^{}
 z_{3,3}^{}z_{1,1}^{}=z_{2,2}^{}z_{1,1}^{}
  z_{3,3}^{}=z_{2,2}^{}z_{2,2}^{*}=q_{2,2}^{}.
\]
Since $q_{1,1}^{}=1$, it follows that $p_{1,1}^{}=q_{2,2}^{}$ is a projection.
From (\ref{eq:comm-unitary1a}), it now follows that $p_{2,2}^{}$ and
$q_{2,1}^{}$ are projections. Consequently $p_{3,2}^{}=1-p_{2,2}^{}$ is  a
projection. Since $q_{2,1}^{}$ is a projection, $p_{2,1}^{}$ also is, and hence
so is $p_{3,1}=1-p_{2,1}^{}-p_{1,1}^{}$. Thus $p_{i,j}^{}$ is a projection for
all $i\geq j$. Hence $q_{i,j}^{}$ also is a projection for all $i\geq j$.

Note that
\[
p_{3,1}^{}=z_{3,1}^{*}z_{3,1}^{}=z_{1,2}^{}z_{2,3}^{}z_{3,1}^{}=
z_{1,2}^{}z_{3,1}^{}z_{2,3}^{}=z_{3,1}^{}z_{1,2}^{}z_{2,3}^{}=
z_{3,1}^{}z_{3,1}^{*}=q_{3,1}^{}.
\]
From the fact that $p_{1,1}^{}=q_{2,2}^{}$ and from the equalities
(\ref{eq:comm-unitary1a}), we get
\begin{IEEEeqnarray}{rClrClrCl}
  p_{3,1}^{} &=& q_{3,1}^{},\quad
  & p_{3,2}^{}&=& q_{3,2}^{}+q_{3,1}^{}, \quad
  & p_{3,3}^{} &=& q_{3,3}^{}+q_{3,2}^{}+q_{3,1}^{},\\
  p_{1,1}^{} &=& q_{2,2}^{},
  &p_{2,1}^{} &=& q_{2,1}^{}-q_{3,1}^{},\quad
  & p_{2,2}^{} &=& q_{3,3}^{}.
\end{IEEEeqnarray}
Therefore it is now enough to show that the family
$\{q_{i,j}^{}: 1\leq j \leq i \leq 3\}$ is commuting.

Since the $q_{i,j}^{}$'s are projections for $i\geq j$, $z_{i,j}^{}$ is a
partial isometry for $i\geq j$. Since $q_{1,1}^{}=I$ and $q_{2,2}^{}$ commute,
it follows from the lemma above that $z_{2,2}^{*}z_{1,1}^{}$ is a partial
isometry. Since $z_{2,2}^{*}z_{1,1}^{}=z_{1,1}^{}z_{2,2}^{*}$, one more
application of the lemma tells us that $p_{1,1}^{}$ and $p_{2,2}^{}$ commute,
i.e.\ $q_{2,2}^{}$ commutes with $q_{3,3}^{}$. From the relation
$z_{2,2}^{}z_{3,1}^{}=z_{3,1}^{}z_{2,2}^{}$, it follows that $q_{2,2}^{}$
commutes with $p_{3,1}^{} = q_{3,1}^{}$. Hence $q_{2,2}^{}$ commutes with
$q_{3,2}^{}=1-q_{3,1}^{}-q_{3,3}^{}$. It now follows that the family
$\{q_{i,j}^{}: 1\leq j \leq i \leq 3\}$ is commuting.
\qed

\begin{corollary}\label{crlre:partial-iso-3}
 Any product of the form $ab$, where $a,b\in \{z_{i,j}^{}: i\geq j\}\cup
\{z_{i,j}^{*}:i\geq j\}$, is a partial isometry.
\end{corollary}
\prf
The proof follows from Lemma~\ref{lemma:partial-iso-1} and
Lemma~\ref{lemma:partial-iso-2}.
\qed

\subsection{Irreducible representations}
We now come to the proof of Theorem~\ref{thm:irred-q0-4}. 
We will split it into different cases, 
corresponding to the different words in the Weyl group $\mathscr{S}_{3}$ 
for the type $A_{2}$ case.
The proof is computational, but the details are instructive and
should provide insight that should help in dealing with the higher rank cases.

Let us start with an  irreducible representation
$\pi$ of $A_2(0)$ on a Hilbert space $\mathcal{H}$. We will denote by $P_{i,j}$
and $Q_{i,j}$ the operators $\pi(p_{i,j}^{})$ and $\pi(q_{i,j}^{})$
respectively. It is easy to see that exactly one of the following conditions will hold for the irreducible representation $\pi$:
\begin{enumerate}
  \item 
  $\pi(z_{3,1})=0$, $\pi(z_{3,2})=0$, $\pi(z_{2,1})=0$,
  \item 
  $\pi(z_{3,1})=0$, $\pi(z_{3,2})=0$, $\pi(z_{2,1})\neq 0$,
  \item 
  $\pi(z_{3,1})=0$, $\pi(z_{3,2})\neq 0$, $\pi(z_{2,1})=0$,
  \item 
  $\pi(z_{3,1})=0$, $\pi(z_{3,2})\neq 0$, $\pi(z_{2,1})\neq 0$,
  \item 
  $\pi(z_{3,1})\neq 0$, $\pi(z_{2,1})\pi(z_{3,2})=0$,
  \item 
  $\pi(z_{3,1})\neq 0$, $\pi(z_{2,1})\pi(z_{3,2})\neq 0$.
\end{enumerate}
We will deal with the above cases separately.
As we will see, these cases correspond to the six elements of the Weyl group $\mathscr{S}_{3}$ in the type $A_{2}$ case, expressed as reduced words in terms of the transpositions $s_{1}=(1,2)$ and $s_{2}=(2,3)$.
\subsubsection{Case 1: $\pi(z_{3,1})=0$, $\pi(z_{3,2})=0$, $\pi(z_{2,1})=0$.}
From the commutation relations 
(\ref{eq:comm-z1}--\ref{eq:comm-z8}), it
follows that $\pi(z_{i,j})=0$ for all $i\neq j$. Therefore $\pi(z_{i,i})$ are
unitary with $\pi(z_{1,1})\pi(z_{2,2})=\pi(z_{3,3}^{*})$ and $\pi(z_{1,1})$
commutes with $\pi(z_{2,2})$. Irreducibility of $\pi$ now implies that
$\mathcal{H}$ is one dimensional and there exist $\lambda,\mu\in S^{1}$ such
that $\pi(z_{1,1})=\lambda$, $\pi(z_{2,2})=\mu$ and
$\pi(z_{3,3})=\overline{\lambda\mu}$.  Thus $\pi$ is unitarily equivalent to
$\psi_{(\lambda,\lambda\mu),id}:=\mychi_{(\lambda,\lambda\mu)}$.

\subsubsection{Case 2: $\pi(z_{3,1})=0$, $\pi(z_{3,2})=0$,
        $\pi(z_{2,1})\neq 0$.}
From the relations 
(\ref{eq:comm-z1}--\ref{eq:comm-z8}), it follows that
\begin{enumerate}
  \item
  $\pi(z_{1.3})=\pi(z_{2,3})=0$,
  \item
  $\pi(z_{3,3})$ is a unitary,
  \item
  $\pi(z_{2,1})$ is normal,
  \item
  $\pi(z_{2,1})$ commutes with $\pi(z_{3,3})$ and $\pi(z_{1,2})$.
\end{enumerate}

\begin{proposition}\label{ppsn:case2a}
Let $P$ be a projection such that
\[
P\leq P_{2,1},\qquad
P\pi(z_{2,1})=\pi(z_{2,1})P,\qquad P\pi(z_{3,3})=\pi(z_{3,3})P.
\]
Then $\mathcal{H}_P:=\{\pi(z_{2,2})^{n}\xi: n\in\mathbb{N},\xi\in
P\mathcal{H}\}$ is an invariant subspace for $\pi$.

If $Q$ is another projection orthogonal to $P$  satisfying the above conditions,
then $\mathcal{H}_{P}$ and $\mathcal{H}_{Q}$ are orthogonal.
\end{proposition}
\prf
From the commutation relations, it follows that
\begin{IEEEeqnarray}{rCl}
  \pi(z_{1,1})\left(\pi(z_{2,2})^n\xi\right)   &=&
  \begin{cases}
  \pi(z_{1,1})\xi  & \text{if  } n=0,\\
  \pi(z_{3,3})^*\pi(z_{2,2})^{n-1}\xi & \text{if  }n>0,
  \end{cases} \label{eq:case2inv-1}\\
  \pi(z_{1,2})\left(\pi(z_{2,2})^n\xi\right)  &=&
  \begin{cases}
  \pi(z_{1,2})\xi & \text{if  }n=0,\\
   0 & \text{if  }n>0,
  \end{cases}  \label{eq:case2inv-2}\\
  \pi(z_{2,1})\left(\pi(z_{2,2})^n\xi\right)  &= &
  \begin{cases}
  \pi(z_{2,1})\xi & \text{if  }n=0,\\
  0&\text{if  }n>0,
  \end{cases}  \label{eq:case2inv-3}\\
  \pi(z_{2,2})\left(\pi(z_{2,2})^n\xi\right) &=
   &\pi(z_{2,2})^{n+1}\xi,\label{eq:case2inv-4}\\
  \pi(z_{3,3})\left(\pi(z_{2,2})^n\xi\right)  &=&
   \pi(z_{2,2})^{n}~\pi(z_{3,3})\xi. \label{eq:case2inv-5}
\end{IEEEeqnarray}
Since $P\leq P_{2,1}$, we have $\xi=P_{2,1}\xi=Q_{2,1}\xi$  for all
$\xi\in\mathcal{H}_{P}$. Hence
\[
\pi(z_{1,1})\xi=\pi(z_{1,1})\pi(z_{2,1})\left(\pi(z_{2,1}^{*})\xi\right)=0.
\]
Thus $\mathcal{H}_{P}$ is invariant under $\pi(z_{1,1})$. Since
$z_{1,2}^{*}=z_{2,1}z_{3,3}$, $\pi(z_{1,2})$ commutes with  $P$ and we  have
$\pi(z_{1,2})\xi=\pi(z_{1,2})P\xi=P\pi(z_{1,2})\xi$.  Thus $\mathcal{H}_{P}$ is
invariant under $\pi(z_{1,2})$. Similarly,  since $\pi(z_{2,1})$ and
$\pi(z_{3,3})$ both commute with $P$, we have the  invariance of
$\mathcal{H}_{P}$ under the actions of $\pi(z_{2,1})$ and $\pi(z_{3,3})$.

From the relations~
(\ref{eq:comm-z1}--\ref{eq:comm-z8}), it follows that
$\mathcal{H}_{P}$ is invariant under $\pi$.

For the second part, take $\xi\in P\mathcal{H}$ and $\zeta\in Q\mathcal{H}$.
Then  for $n\in\mathbb{N}$, one has
\begin{IEEEeqnarray*}{rCl}
  \langle \pi(z_{2,2})^{n}\zeta, \pi(z_{2,2})^{n}\xi\rangle &=&
  \langle \pi(z_{2,2}^{*})^{n}\pi(z_{2,2})^{n}\zeta,\xi\rangle\\
  &=& \langle \zeta,\xi\rangle\\
  &=&\langle Q\zeta , P\xi \rangle\\
  &=& 0.
\end{IEEEeqnarray*}
Next, for $m,n\in\mathbb{N}$ with $m>n$, we have
\begin{IEEEeqnarray*}{rCl}
\langle \pi(z_{2,2})^{m}\zeta, \pi(z_{2,2})^{n}\xi\rangle &=&
\langle \zeta, \pi(z_{2,2}^{*})^{m}\pi(z_{2,2})^{n}\xi\rangle\\
  &=& \langle \zeta,\pi(z_{2,2}^{*})^{m-n}\xi\rangle\\
  &=&\langle \zeta, \pi(z_{3,3})^{m-n}\pi(z_{1,1})^{m-n}\xi\rangle\\
  &=&0.
\end{IEEEeqnarray*}
A similar calculation gives $\langle \pi(z_{2,2})^{m}\zeta,
\pi(z_{2,2})^{n}\xi\rangle=0$ for $m<n$. Thus $\mathcal{H}_{P}$ and
$\mathcal{H}_{Q}$ are orthogonal.
\qed

\begin{corollary}\label{crlre:case2b}
There exists a $\lambda\in S^{1}$ such that
$\sigma(\pi(z_{2,1}))=\{0,\lambda\}$.
\end{corollary}
\prf
Let $E$ and $F$ be two disjoint closed subsets of the spectrum
$\sigma(\pi(z_{2,1}))$ such that $0\notin E$, $0\notin F$ and let $P$ and $Q$ be
the spectral projections of $\pi(z_{2,1})$ corresponding to the subsets $E$ and
$F$ respectively. Then it follows from the previous proposition that
$\mathcal{H}_{P}$ and $\mathcal{H}_{Q}$ are two invariant subspaces for $\pi$
that are mutually orthogonal. By irreducibility of $\pi$, it follows that
$\sigma(\pi(z_{2,1}))=\{0,\lambda\}$ for some $\lambda\neq 0$ in $\mathbb{C}$.
Since $\pi(z_{2,1})$ is a normal partial isometry, it follows that $\lambda\in
S^{1}$.
\qed

\begin{corollary}\label{crlre:case2c}
There exists a $\mu\in S^{1}$ such that
$\sigma(\pi(z_{3,3}))=\{\mu\}$.
\end{corollary}
\prf
From the commutation relations, it follows that $\pi(z_{3,3})$ is a unitary and
commutes with $P_{2,1}$. Let $P$ and $Q$ be two spectral projections
corresponding to two disjoint closed subsets of the spectrum of the restriction
of $\pi(z_{3,3})$ to $P_{2,1}\mathcal{H}$. Then the above proposition tells us
that $\mathcal{H}_{P}$ and $\mathcal{H}_{Q}$ are two invariant subspaces for
$\pi$ that are mutually orthogonal. Using irreducibility of $\pi$, we deduce
that $\sigma(\pi(z_{3,3})|_{P_{2,1}\mathcal{H}})=\{\mu\}$ where $\mu\in S^{1}$.
Hence for $\xi\in P_{2,1}\mathcal{H}$, we have
\[
\pi(z_{3,3}^{*})(\pi(z_{2,2}^{n})\xi)=(\pi(z_{2,2}^{n})\pi(z_{3,3}^{*})\xi=
\bar{\mu}\pi(z_{2,2}^{n})\xi.
\]
Since $\mathcal{H}_{P_{2,1}}=\mathcal{H}$, we have the result.
\qed

\begin{proposition}\label{ppsn:case2d}
The projection $P_{2,1}$ is of rank one and $\pi$ is unitarily equivalent to the
representation $\psi_{(\bar{\lambda}\bar{\mu},\bar{\mu}),s_{1}}$.
\end{proposition}
\prf
Take a unit vector $\xi\in P_{2,1}\mathcal{H}$. From the previous two
corollaries, it follows that $\pi(z_{3,3})\xi=\mu\xi$ and
$\pi(z_{2,1})\xi=\lambda\xi$. Therefore
$\pi(z_{1,2})\xi=\pi(z_{3,3}^{*}z_{2,1}^{*})\xi=\bar{\lambda}\bar{\mu}\xi$.
Let
\[
\mathcal{H}_{\xi}=\text{span\,}\{\pi(z_{2,2})^{n}\xi: n\in\mathbb{N}\}.
\]
From equations (\ref{eq:case2inv-1}--\ref{eq:case2inv-5}), it now follows that
\begin{IEEEeqnarray}{rCl}
  \pi(z_{1,1})\left(\pi(z_{2,2})^n\xi\right)   &=&
  \begin{cases}
  0  &\text{if  } n=0,\\
  \bar{\mu}\pi(z_{2,2})^{n-1}\xi & \text{if  }n>0,
  \end{cases} \label{eq:case2inv-6}\\
  \pi(z_{1,2})\left(\pi(z_{2,2})^n\xi\right)  &=&
  \begin{cases}
  \bar{\lambda}\bar{\mu}\xi & \text{if  }n=0,\\
   0 & \text{if  }n>0,
  \end{cases}  \label{eq:case2inv-7}\\
  \pi(z_{2,1})\left(\pi(z_{2,2})^n\xi\right)  &= &
  \begin{cases}
  \lambda\xi & \text{if  }n=0,\\
  0&\text{if  }n>0,
  \end{cases}  \label{eq:case2inv-8}\\
  \pi(z_{2,2})\left(\pi(z_{2,2})^n\xi\right) &= &
  \pi(z_{2,2})^{n+1}\xi,\label{eq:case2inv-9}\\
  \pi(z_{3,3})\left(\pi(z_{2,2})^n\xi\right)  &=&
  \mu\pi(z_{2,2})^{n}\xi. \label{eq:case2inv-10}
\end{IEEEeqnarray}
Thus $\mathcal{H}_{\xi}$ is an invariant subspace of $\mathcal{H}$ and hence by
irreducibility of $\pi$, one has $\mathcal{H}_{\xi}=\mathcal{H}$.

The map
\[
U: e_{n}\mapsto \lambda^{-n} \pi(z_{2,2})^n\xi
\]
from $\ell^{2}(\mathbb{N})$ to $\mathcal{H}$ now extends to a unitary and gives
us the required unitary equivalence.
\qed

\subsubsection{Case 3: $\pi(z_{3,1})=0$, $\pi(z_{3,2})\neq 0$,
   $\pi(z_{2,1})=0$.}
Analogous to the results in the previous case, here we have the following
results. The proofs are similar.

\begin{proposition}\label{ppsn:case3a}
Let $P$ be a projection such that
\[
P\leq P_{3,2},\qquad P\pi(z_{3,2})=\pi(z_{3,2})P, \qquad
 P\pi(z_{1,1})=\pi(z_{1,1})P.
\]
Then $\mathcal{H}_P:=\{\pi(z_{3,3})^{n}\xi: n\in\mathbb{N},\xi\in
P\mathcal{H}\}$ is an invariant subspace for $\pi$.

If $Q$ is another projection orthogonal to $P$  satisfying the conditions above,
then $\mathcal{H}_{P}$ and $\mathcal{H}_{Q}$ are orthogonal.
\end{proposition}

\begin{corollary}\label{crlre:case3b}
There exists a $\lambda\in S^{1}$ such that
$\sigma(\pi(z_{3,2}))=\{0,\lambda\}$.
\end{corollary}

\begin{corollary}\label{crlre:case3c}
There exists a $\mu\in S^{1}$ such that
$\sigma(\pi(z_{1,1}))=\{\mu\}$.
\end{corollary}

\begin{proposition}\label{ppsn:case3d}
The projection $P_{3,2}$ is of rank one and $\pi$ is unitarily equivalent to
the representation $\psi_{(\mu,\bar{\lambda}),s_{2}}$.
\end{proposition}

\subsubsection{Case 4: $\pi(z_{3,1})=0$, $\pi(z_{3,2})\neq 0$,
$\pi(z_{2,1})\neq 0$.}
Let us first prove the following lemma.

\begin{lemma}\label{lemma:case4a}
If $\pi(z_{3,1})=0$ and $\pi(z_{2,1})\pi(z_{3,2})=0$, then either
$\pi(z_{2,1})=0$ or $\pi(z_{3,2})=0$.
\end{lemma}
\prf
Let us assume that $\pi(z_{3,1})=0$, $\pi(z_{2,1})\pi(z_{3,2})=0$ and
$\pi(z_{2,1})\neq 0$. We will show that $\pi(z_{3,2})=0$. Since
$\pi(z_{3,1})=0$, it follows from 
(\ref{eq:comm-z1}--\ref{eq:comm-z3})
that $\pi(z_{1,1})\pi(z_{3,2})=\pi(z_{3,2})\pi(z_{1,1})$.
Since $z_{2,3}$ commutes with $z_{3,2}$, we get
\[
\pi(z_{3,2}^{})\pi(z_{3,2}^{*})=\pi(z_{3,2}z_{1,1}z_{2,3})=
\pi(z_{1,1}z_{2,3}z_{3,2})=\pi(z_{3,2}^{*}z_{3,2}^{}),
\]
i.e.\ $\pi(z_{3,2}^{})$ is normal. Since $\pi(z_{2,1})^*$ commutes with
$\pi(z_{3,2})$, it follows that $\pi(z_{2,1})$ commutes with $\pi(z_{3,2})$.

Let $\mathcal{H}_{0}$ be the closed linear span of
$\{\pi(z_{2,2}^{k})\pi(z_{2,1}^{})\xi: k\in\mathbb{N}, \xi\in \mathcal{H}\}$.
Using the commutation relations
(\ref{eq:comm-z1}--\ref{eq:comm-z8}) and
the fact that $\pi(z_{3,1}^{})=0$, 
one can compute the quantities 
$\pi(z_{i,j})\left(\pi(z_{2,2}^{k})\pi(z_{2,1}^{})\xi\right)$
from which it follows that 
\begin{enumerate}
  \item
  $\mathcal{H}_0$ is invariant subspace for $\pi$, and
  \item
  $\pi(z_{3,2})|_{\mathcal{H}_0}=0$
\end{enumerate}
Since $\pi$ is irreducible, we
have $\mathcal{H}=\mathcal{H}_0$ and $\pi(z_{3,2})=0$.
\qed

\begin{proposition}\label{ppsn:case4b}
  Let $P$ be a projection such that
  \[
  P\leq P_{1,3},\qquad P\pi(z_{2,1})=\pi(z_{2,1})P,\qquad
  P\pi(z_{3,2})=\pi(z_{3,2})P.
  \]
  Then the subspace $\mathcal{H}_P:=\{\pi(z_{1,1}^*)^{m}\pi(z_{3,3})^{n}\xi:
  m,n\in\mathbb{N},\xi\in P\mathcal{H}\}$  is invariant for $\pi$.

 If $Q$ is another projection orthogonal to $P$  satisfying the same conditions
 above, then $\mathcal{H}_{P}$ and $\mathcal{H}_{Q}$ are orthogonal.
\end{proposition}
\prf
Let $\xi\in P\mathcal{H}$. Then we have
\begin{equation}\label{eq:temp1}
\pi(z_{1,1})\xi=0=\pi(z_{1,2})\xi=\pi(z_{2,2})\xi=\pi(z_{2,2}^*)\xi=
\pi(z_{2,3}^*)\xi=\pi(z_{3,3}^*)\xi.
\end{equation}
Therefore using the given conditions, one can show that
$\mathcal{H}_P$ is an invariant subspace for $\pi$. Irreducibility
of $\pi$ now gives $\mathcal{H}=\mathcal{H}_P$.

For the second part, Let us compute the inner product $\langle
\pi(z_{1,1}^*)^{m}\pi(z_{3,3})^{n}\xi,
\pi(z_{1,1}^*)^{m'}\pi(z_{3,3})^{n'}\xi'\rangle$ where $\xi\in P\mathcal{H}$ and
$\xi'\in Q\mathcal{H}$. In the case $m\neq m'$ or $n\neq n'$, using  the
relations $z_{1,1}^{}z_{1,1}^{*}=1=z_{3,3}^{*}z_{3,3}^{}$,
$z_{1,1}^{}z_{3,3}^{}=z_{3,3}^{}z_{1,1}^{}$ and
$z_{1,1}^{}z_{3,3}^{*}=z_{3,3}^{*}z_{1,1}^{}$, one arrives at a vector of the
form $\pi(a)\zeta$  where $a$ is $z_{1,1}^{*}$ or $z_{3,3}^{*}$ and $\zeta$ is
either $\xi$ or $\xi'$, and by (\ref{eq:temp1}), $\pi(a)\zeta=0$. If $m=m'$ and
$n=n'$,  the above inner product reduces to $\langle \xi,\xi'\rangle$, which is zero.
\qed

Proof of the next two corollaries are similar to the earlier cases and hence
omitted.
\begin{corollary}\label{crlre:case4c}
The operator $\pi(z_{2,1}^{})$ is normal, with $P_{1,3}\leq P_{2,1}$ and there
exists a $\lambda\in S^{1}$ such that the restriction of $\pi(z_{2,1})$ to
$P_{1,3}\mathcal{H}$ is $\lambda I$.
\end{corollary}

\begin{corollary}\label{crlre:case4d}
The operator $\pi(z_{3,2}^{})$ is normal, with $P_{1,3}\leq P_{3,2}$ and there
exists a $\mu\in S^{1}$ such that the restriction of $\pi(z_{3,2})$ to
$P_{1,3}\mathcal{H}$ is $\mu I$.
\end{corollary}

\begin{proposition}\label{ppsn:case4e}
The representation $\pi$ is unitarily equivalent to
$\psi_{(\bar{\lambda}\bar{\mu},\bar{\mu}),s_{1}s_{2}}$.
\end{proposition}
\prf
Take a unit vector $\xi\in P_{1,3}\mathcal{H}$. From the previous two
corollaries, it follows that $\pi(z_{3,2}^{})\xi=\mu\xi$ and
$\pi(z_{2,1}^{}z_{3,2}^{})\xi=\lambda\mu\xi$.
Let
\[
\mathcal{H}_{\xi}=\{\pi(z_{1,1}^*)^{m}\pi(z_{3,3})^{n}\xi: m,n\in\mathbb{N}\}.
\]
From the computations in the proof of Proposition~\ref{ppsn:case4b},
it now follows that
\begin{IEEEeqnarray*}{rCl}
    \pi(z_{1,1})\left(\pi(z_{1,1}^*)^{m}\pi(z_{3,3})^{n}\xi\right)&=&
    \begin{cases}
    0&\text{if  } m=0,\\
    \pi(z_{1,1}^*)^{m-1}\pi(z_{3,3})^{n}\xi&\text{if  } m>0,
    \end{cases}\\
    &&\\
  \pi(z_{1,2})\left(\pi(z_{1,1}^*)^{m}\pi(z_{3,3})^{n}\xi\right)  &=&
  \begin{cases}
     0 & \text{if  } m=0, n=0,\\
  \pi(z_{3,3})^{n-1} \left(\bar{\mu}\xi\right) &\text{if  } m=0, n>0,\\
  0&\text{if  }m>0,
  \end{cases}\\
    &&\\
  \pi(z_{1,3})\left(\pi(z_{1,1}^*)^{m}\pi(z_{3,3})^{n}\xi\right) &=&
  \begin{cases}
  \bar{\mu}\bar{\lambda}\xi &\text{if  } m=0, n=0,\\
  0&\text{if  }m=0,n>0,\\
  0&\text{if  }m>0,
  \end{cases} \\
    &&\\
  \pi(z_{2,1})\left(\pi(z_{1,1}^*)^{m}\pi(z_{3,3})^{n}\xi\right) &=&
  \begin{cases}
  \pi(z_{3,3})^{n} \left(\lambda\xi\right) &\text{if  } m=0,\\
  0&\text{if  }m>0,
  \end{cases}\\
    &&\\
  \pi(z_{2,2})\left(\pi(z_{1,1}^*)^{m}\pi(z_{3,3})^{n}\xi\right) &=&
   \begin{cases}
  0&\text{if  } n=0,\\
  \pi(z_{1,1}^*)^{m+1}\pi(z_{3,3})^{n-1} \xi &\text{if  } n>0,
  \end{cases}\\
    &&\\
  \pi(z_{2,3})\left(\pi(z_{1,1}^*)^{m}\pi(z_{3,3})^{n}\xi\right) &=&
  \begin{cases}
  \pi(z_{1,1}^*)^{m+1} \left(\bar{\mu}\xi\right) &\text{if  } n=0,\\
  0&\text{if  }n>0,
  \end{cases}\\
    &&\\
  \pi(z_{3,2})\left(\pi(z_{1,1}^*)^{m}\pi(z_{3,3})^{n}\xi\right) &=&
  \begin{cases}
  \pi(z_{1,1}^*)^{m} \left(\bar{\mu}\xi\right) &\text{if  } n=0,\\
  0&\text{if  }n>0,
  \end{cases}\\
    &&\\
  \pi(z_{3,3})\left(\pi(z_{1,1}^*)^{m}\pi(z_{3,3})^{n}\xi\right)
  &=&\pi(z_{1,1}^*)^{m}\pi(z_{3,3})^{n+1} \xi.
\end{IEEEeqnarray*}
Thus $\mathcal{H}_{\xi}$ is an invariant subspace of $\mathcal{H}$ and hence by
irreducibility of $\pi$, one has $\mathcal{H}_{\xi}=\mathcal{H}$.

The map
\[
U: e_{m,n}^{}\mapsto \mu^{-m-n}\lambda^{-m} \pi(z_{1,1}^*)^{m}\pi(z_{3,3})^n \xi
\]
from $\ell^{2}(\mathbb{N}^{2})$ to $\mathcal{H}$ gives us the required unitary.
\qed

\subsubsection{Case 5: $\pi(z_{3,1})\neq 0$, $\pi(z_{2,1})\pi(z_{3,2})=0$.}
Note that in both Case 5 and Case 6, one has
$\pi(z_{3,1}^{*})=\pi(z_{1,2})\pi(z_{2,3})\neq 0$. This implies that
$\pi(z_{2,3})\neq 0$, which in turn implies that $\pi(z_{3,2})\neq 0$ because
$\pi(z_{2,3}^{*})=\pi(z_{1,1})\pi(z_{3,2})$.

\begin{proposition}\label{ppsn:case5b}
  Let $P$ be a projection such that $P\leq P_{3,1}$ and
  $P\pi(z_{1,2})=\pi(z_{1,2})P$,  $P\pi(z_{2,3})=\pi(z_{2,3})P$.  Let
  $\mathcal{H}_P:=\{\pi(z_{1,1}^*)^{n}\pi(z_{3,3})^{m}\xi:n,m\geq 0,\xi\in
  P\mathcal{H}\}$. Then $\mathcal{H}_P$ is an invariant subspace for $\pi$.

  If $Q$ is another projection orthogonal to $P$  such that
  $Q\pi(z_{1,2})=\pi(z_{1,2})Q$ and  $Q\pi(z_{2,3})=\pi(z_{2,3})Q$, then
  $\mathcal{H}_{P}$ and $\mathcal{H}_{Q}$ are orthogonal.
\end{proposition}

\begin{corollary}\label{crlre:case5c}
The operator $\pi(z_{1,2}^{})$ is normal, with $P_{3,1}\leq P_{1,2}$ and there
exists a $\lambda\in S^{1}$ such that the restriction of $\pi(z_{1,2}^{})$ to
$P_{3,1}\mathcal{H}$ is $\lambda I$.
\end{corollary}

\begin{corollary}\label{crlre:case5d}
The operator $\pi(z_{2,3}^{})$ is normal, with $P_{3,1}\leq P_{2,3}$ and there
exists a $\mu\in S^{1}$ such that the restriction of $\pi(z_{2,3}^{})$ to
$P_{3,1}\mathcal{H}$ is $\mu I$.
\end{corollary}

\begin{proposition}\label{ppsn:case5e}
The representation $\pi$ is unitarily equivalent to
$\psi_{(\lambda,\lambda\bar{\mu}),s_{2}s_{1}}$.
\end{proposition}

\subsubsection{Case 6: $\pi(z_{3,1})\neq 0$, $\pi(z_{2,1})\pi(z_{3,2})\neq 0$.}

\begin{lemma}\label{lemma:case6a}
If $\pi(z_{3,1})\neq 0$ and $\pi(z_{2,1})\pi(z_{3,2})\neq 0$,
then $\pi(z_{3,1}^{}z_{2,1}^{}z_{3,2}^{})\neq 0$.
\end{lemma}
\prf
Let us assume that $\pi(z_{3,1})\neq 0$ and
$\pi(z_{3,1}^{}z_{2,1}^{}z_{3,2}^{})= 0$. We will show that this implies
$\pi(z_{2,1})\pi(z_{3,2})= 0$.

Let $\mathcal{H}_{3,1}$ be the closed linear span of
$\left\{\pi(z_{1,1}^{*})^{m}\pi(z_{3,3}^{n})\xi:m,n\in\mathbb{N},\xi\in
P_{3,1}\mathcal{H}\right\}$. We will show that $\mathcal{H}=\mathcal{H}_{3,1}$
and $\pi(z_{2,1})\pi(z_{3,2})=0$.

Let $\xi\in P_{3,1}\mathcal{H}$. Then one has
$\xi=P_{3,1}\xi=\pi(z_{3,1}^{*})\xi'$ for some $\xi'\in\mathcal{H}$. Since
$z_{3,1}^{}$ is normal, one also has $\xi=\pi(z_{3,1}^{})\xi''$ for some
$\xi''\in\mathcal{H}$. Using this observation, the relation
$z_{1,3}^{*}=z_{2,1}^{}z_{3,2}^{}$ and the condition
$\pi(z_{3,1}^{}z_{2,1}^{}z_{3,2}^{})= 0$, we now get
\begin{IEEEeqnarray}{rClrCl}
\pi(z_{1,3}^{})\xi &=& 0 &\qquad \pi(z_{1,3}^{*})\xi &= & 0,\\
\pi(z_{1,1}^{})\xi &=& 0 &\qquad \pi(z_{2,1}^{})\xi &= & 0,\\
\pi(z_{3,2}^{*})\xi &=& 0 &\qquad \pi(z_{3,3}^{*})\xi &= & 0.
\end{IEEEeqnarray}
Since $z_{1,2}^{}$ and $z_{2,3}^{}$ commute with $z_{3,1}^{*}$ and $z_{3,1}^{}$
is normal, it also follows that
\begin{IEEEeqnarray}{rClCrCl}
\pi(z_{1,2}^{})\xi &\in & \mathcal{H}_{3,1} &&
\qquad \pi(z_{1,2}^{*})\xi &\in & \mathcal{H}_{3,1},\\
\pi(z_{2,3}^{})\xi &\in & \mathcal{H}_{3,1} &&
\qquad \pi(z_{2,3}^{*})\xi &\in & \mathcal{H}_{3,1},\\
\pi(z_{2,3}^{*})\pi(z_{1,1}^{*})^{}\xi 
&=&
 \pi(z_{1,1}^{*})^{}\pi(z_{2,3}^{*})\xi,&& \qquad \xi &\in& P_{3,1}\mathcal{H}.
\end{IEEEeqnarray}
Using the above, one can now verify that
\begin{enumerate}
  \item
  $\mathcal{H}_{3,1}$ is invariant subspace for $\pi$, and
  \item
  $\pi(z_{1,3})|_{\mathcal{H}_{3,1}}=0$
\end{enumerate}
By irreducibility of $\pi$, $\mathcal{H}=\mathcal{H}_{3,1}$ 
and $\pi(z_{1,3})=0$. Therefore
$\pi(z_{2,1})\pi(z_{3,2})=\pi(z_{1,3}^{*})=0$.
\qed

Observe that the relations 
(\ref{eq:comm-z1}--\ref{eq:comm-z8})
imply
that $z_{1,3}^{}$ and $z_{3,1}^{}$ are both normal, i.e.\ $P_{1,3}=Q_{1,3}$ and
$P_{3,1}=Q_{3,1}$. Using the relations (\ref{eq:comm-z1}--\ref{eq:comm-z8}) again,
we now
conclude that $z_{3,1}^{}z_{2,1}^{}z_{3,2}^{}$ is a normal partial isometry,
with initial and final projection given by $P_{1,3}P_{3,1}$ as the following
computations show
\begin{IEEEeqnarray*}{rCl}
  \left(z_{3,1}^{}z_{2,1}^{}z_{3,2}^{}\right)
        \left(z_{3,1}^{}z_{2,1}^{}z_{3,2}^{}\right)^{*} &=&
    \left(z_{3,1}^{}z_{1,3}^{*}\right)  \left(z_{3,1}^{}z_{1,3}^{*}\right)^{*}\\
  &=& z_{3,1}^{}z_{1,3}^{*}z_{1,3}^{}z_{3,1}^{*}\\
  &=& z_{1,3}^{*}z_{1,3}^{}z_{3,1}^{}z_{3,1}^{*}\\
  &=& z_{1,3}^{*}z_{1,3}^{}z_{3,1}^{*}z_{3,1}^{}
\end{IEEEeqnarray*}
and similarly,
\begin{IEEEeqnarray*}{rCl}
  \left(z_{3,1}^{}z_{2,1}^{}z_{3,2}^{}\right)^{*}
     \left(z_{3,1}^{}z_{2,1}^{}z_{3,2}^{}\right)&=&
  \left(z_{3,1}^{}z_{1,3}^{*}\right)^{*}  \left(z_{3,1}^{}z_{1,3}^{*}\right)  \\
  &=& z_{1,3}^{}z_{3,1}^{*}z_{3,1}^{}z_{1,3}^{*}\\
  &=& z_{1,3}^{}z_{1,3}^{*}z_{3,1}^{*}z_{3,1}^{}\\
  &=& z_{1,3}^{*}z_{1,3}^{}z_{3,1}^{*}z_{3,1}^{}.
\end{IEEEeqnarray*}

\begin{proposition}\label{ppsn:case6b}
  Let $P$ be a projection such that
  \[
  P\leq P_{1,3}P_{3,1},\qquad P\pi(z_{3,1}^{})=\pi(z_{3,1}^{})P,\qquad
  P\pi(z_{1,3}^{})=\pi(z_{1,3}^{})P.
  \]
  Then
  $\mathcal{H}_P:=\left\{\pi(z_{3,3}^{m})\pi(z_{2,3}^{k})\pi(z_{3,2}^{n})\xi:
      k,m,n\geq 0,\xi\in P\mathcal{H}\right\}$
      is an invariant subspace for $\pi$.

  If $Q$ is another projection orthogonal to $P$ satisfying the same conditions
  above, then $\mathcal{H}_{P}$ and $\mathcal{H}_{Q}$ are orthogonal.
\end{proposition}
\prf
Since $P\leq P_{1,3}P_{3,1}$, one has $\xi=P_{1,3}P_{3,1}\xi$ for $\xi\in
P\mathcal{H}$. Therefore
\begin{IEEEeqnarray}{rClrClrCl}
  \pi(z_{1,1}^{})\xi &=& 0,\qquad &\pi(z_{1,2}^{})\xi
  &=& 0,\qquad & \pi(z_{2,1}^{})\xi &=& 0,\label{eq:repn-1}\\
  \pi(z_{2,3}^{*})\xi &=& 0,\qquad & \pi(z_{3,2}^{*})\xi
  &=& 0,\qquad & \pi(z_{3,3}^{*})\xi &=& 0.\label{eq:repn-2}
\end{IEEEeqnarray}
Consequently, we also have
\begin{equation}\label{eq:repn-3}
   \pi(z_{2,2}^{})\xi=\pi(z_{3,3}^{*})\pi(z_{1,1}^{*})\xi=
   \pi(z_{1,1}^{*})\pi(z_{3,3}^{*})\xi=0.
\end{equation}

From (\ref{eq:repn-1}-\ref{eq:repn-2}) and 
(\ref{eq:comm-z1}--\ref{eq:comm-z3}),
it follows that for all $n\in\mathbb{N}$.
\begin{equation}\label{eq:repn-4}
\pi(z_{2,1})\pi(z_{2,3}^{n})\xi  = 0, \qquad
\pi(z_{1,2})\pi(z_{3,2}^{n})\xi  = 0.
\end{equation}

Since $P$ commutes with $\pi(z_{1,3}^{})$ and  $\pi(z_{3,1}^{})$, we have
\begin{equation}\label{eq:repn-5}
\xi\in P\mathcal{H} \Longrightarrow \pi(z_{1,3}^{})\xi,
\pi(z_{1,3}^{*})\xi, \pi(z_{3,1}^{})\xi, \pi(z_{3,1}^{*})\xi \in P\mathcal{H}.
\end{equation}

Let $\xi\in P\mathcal{H}$.

\begin{IEEEeqnarray*}{rCl}
   \pi(z_{1,1}^{})
        \left(\pi(z_{3,3}^{m})\pi(z_{2,3}^{k})\pi(z_{3,2}^{n})\xi\right) 
    &=&
    \begin{cases}
    0 & \text{if  } k=0,\\
    \pi(z_{3,3}^{m})\pi(z_{2,3}^{k-1})\Bigl(\pi(z_{3,2}^{*})
    \pi(z_{3,2}^{n})\xi\Bigr)&    \text{if  } k>0,
    \end{cases}\yesnumber\label{eq:action-6a-11}\\
    &&\\
    \pi(z_{1,2})
          \left(\pi(z_{3,3}^{m})\pi(z_{2,3}^{k})\pi(z_{3,2}^{n})\xi\right) 
    &=&
     \begin{cases}
    0 &\text{if  }m=0,k=0,\\
     \pi(z_{2,3}^{k-1})\Bigl(\pi(z_{3,1}^{*})
     \pi(z_{3,2}^{n})\xi\Bigr)&\text{if  }m=0,k>0,\\
     \pi(z_{3,3}^{m-1})\pi(z_{3,2}^{n+1})\Bigl(\pi(z_{1,3}^{})\xi\Bigr)
     &\text{if  }m>0,k=0,\\
     \pi(z_{3,3}^{m})\pi(z_{2,3}^{k-1})\Bigl(\pi(z_{3,1}^{*})
     \pi(z_{3,2}^{n})\xi\Bigr) & \text{if  }m>0,k>0,
       \end{cases}\\
       &&\yesnumber\label{eq:action-6a-12}\\
  \pi(z_{1,3}^{})
        \left(\pi(z_{3,3}^{m})\pi(z_{2,3}^{k})\pi(z_{3,2}^{n})\xi\right) 
  &=&
  \begin{cases}
    \pi(z_{3,2}^{n}) \Bigl(\pi(z_{1,3}^{})\xi\Bigr) &\text{if  } m=0, k=0,\\
    0&\text{if  } m>0\text{ or } k>0,
  \end{cases}\yesnumber\label{eq:action-6a-13}\\
  &&\\
  \pi(z_{2,1}^{})
       \left(\pi(z_{3,3}^{m})\pi(z_{2,3}^{k})\pi(z_{3,2}^{n})\xi\right) 
  &=&
  \begin{cases}
    0 & \text{if  }m=0, k=0,n=0,\\
  \pi(z_{3,2}^{n-1})\Bigl(\pi(z_{1,3}^{*})\xi\Bigr)&\text{if  }m=0, k=0,n>0,\\
  0&\text{if  }m=0, k>0,\\
      \pi(z_{3,3}^{m-1})\pi(z_{2,3}^{k+1})\Bigl(\pi(z_{3,1}^{})\xi\Bigr)
      &\text{if  }m>0,n=0,\\
  \pi(z_{3,3}^{m})\pi(z_{3,2}^{n-1})
  \Bigl(\pi(z_{1,3}^{*})\xi\Bigr) &\text{if  }m>0,n>0, k=0,\\
  0&\text{if  }m>0,n>0, k>0,
    \end{cases}\\
    &&\yesnumber\label{eq:action-6a-21}\\
  \pi(z_{2,2}^{})
  \left(\pi(z_{3,3}^{m})\pi(z_{2,3}^{k})\pi(z_{3,2}^{n})\xi\right)
  &=&
  \begin{cases}
    0 &\text{if  }m=0,\\
    \pi(z_{3,3}^{m-1})\pi(z_{2,3}^{k+1})\pi(z_{3,2}^{n+1})\xi & \text{if  } m>0,
  \end{cases}\yesnumber\label{eq:action-6a-22}\\
  &&\\
  \pi(z_{2,3}^{})
   \left(\pi(z_{3,3}^{m})\pi(z_{2,3}^{k})\pi(z_{3,2}^{n})\xi\right) 
  &=&
  \begin{cases}
  \pi(z_{2,3}^{k+1})\pi(z_{3,2}^{n}) \xi &\text{if  }m=0,\\
  0&\text{if  }m>0,
    \end{cases}\yesnumber\label{eq:action-6a-23}\\
  &&\\
  \pi(z_{3,1}^{})
  \left(\pi(z_{3,3}^{m})\pi(z_{2,3}^{k})\pi(z_{3,2}^{n})\xi\right) 
  &=&
  \begin{cases}
    \pi(z_{2,3}^{k}) \Bigl(\pi(z_{3,1}^{})\xi\Bigr) &\text{if  } m=0,n=0,\\
    0&\text{if  } m+n>0,
  \end{cases}\yesnumber\label{eq:action-6a-31}\\
  &&\\
  \pi(z_{3,2}^{})
  \left(\pi(z_{3,3}^{m})\pi(z_{2,3}^{k})\pi(z_{3,2}^{n})\xi\right)
  &=&
  \begin{cases}
    \pi(z_{2,3}^{k})\pi(z_{3,2}^{n+1}) \xi &\text{if  }m=0,\\
    0&\text{if  }m>0,
  \end{cases}\yesnumber\label{eq:action-6a-32}\\
  &&\\
  \pi(z_{3,3}^{})
   \left(\pi(z_{3,3}^{m})\pi(z_{2,3}^{k})\pi(z_{3,2}^{n})\xi\right) 
  &=& \pi(z_{3,3}^{m+1})\pi(z_{2,3}^{k})\pi(z_{3,2}^{n}) \xi .
    \yesnumber\label{eq:action-6a-33}
\end{IEEEeqnarray*}

Now observe that
\begin{IEEEeqnarray*}{rCl}
   \pi(z_{3,2}^{*})\pi(z_{3,2}^{n})\xi &=&
     \begin{cases}
     \pi(z_{3,2}^{*})\xi=\pi(z_{3,2}^{*})\pi(z_{3,1}^{*})\pi(z_{3,1}^{})\xi=0
     &\text{if }n=0,\\
     (I-\pi(z_{2,2}^{*})\pi(z_{2,2}^{}))\pi(z_{3,2}^{n-1})\xi=
     \pi(z_{3,2}^{n-1})\xi &\text{if }n>0,
     \end{cases}\\
    \pi(z_{3,1}^*)\pi(z_{3,2}^{n})\xi &=&
          \begin{cases}
        \pi(z_{3,1}^{*})\xi &\text{if }n=0,\\
        \pi(z_{1,2}^{})\pi(z_{2,3}^{})\pi(z_{3,2}^{n})\xi
        =\pi(z_{1,2}^{})\pi(z_{3,2}^{n})\pi(z_{2,3}^{})\xi =0 &\text{if }n>0,
     \end{cases}\\
     \pi(z_{2,1})\pi(z_{3,2}^{n})\xi &=&
          \begin{cases}
        \pi(z_{2,1}^{})\xi=0 &\text{if }n=0,\\
        \pi(z_{1,3}^{*})\pi(z_{3,2}^{n-1})\xi =
        \pi(z_{3,2}^{n-1})\pi(z_{1,3}^{*})\xi  &\text{if }n>0.
           \end{cases}
\end{IEEEeqnarray*}
Therefore it follows that
$\pi(z_{i,j}^{})\Bigl(\pi(z_{3,3}^{m})
\pi(z_{2,3}^{k})\pi(z_{3,2}^{n})\xi\Bigr)\in \mathcal{H}_{P}$ for all $i,j$.
Thus $\mathcal{H}_{P}$ is  an invariant subspace for $\pi$.

For the second part, take $\xi\in P\mathcal{H}$  and $\xi'\in Q\mathcal{H}$. Let
\begin{IEEEeqnarray*}{rCl}
  \xi_{k,m,n}^{} &=& \pi(z_{3,3}^{m}) \pi(z_{2,3}^{k}) \pi(z_{3,2}^{n}) \xi,
  \quad k,m,n\in\mathbb{N},\\
  {\xi'}_{k,m,n}^{} &=& \pi(z_{3,3}^{m}) \pi(z_{2,3}^{k}) \pi(z_{3,2}^{n}) \xi',
  \quad k,m,n\in\mathbb{N}.
\end{IEEEeqnarray*}

Since $z_{3,3}^{*}z_{3,3}^{}=1$, we have
\begin{IEEEeqnarray*}{rCl}
  \langle \xi_{k,m,n}^{},  {\xi'}_{k',m',n'}^{} \rangle
  &=&
  \begin{cases}
    \langle   \pi(z_{3,3}^{*})^{m'-m}\pi(z_{2,3}^{k}) \pi(z_{3,2}^{n}) \xi
    \;, \;  \pi(z_{2,3}^{k'}) \pi(z_{3,2}^{n'}) \xi' \rangle &\text{if }m<m',\\
    \langle   \pi(z_{2,3}^{k}) \pi(z_{3,2}^{n}) \xi \;, \;
    \pi(z_{3,3}^{*})^{m-m'} \pi(z_{2,3}^{k'}) \pi(z_{3,2}^{n'}) \xi' \rangle
           &\text{if }m>m',\\
    \langle  \pi(z_{2,3}^{k}) \pi(z_{3,2}^{n}) \xi \;, \;
    \pi(z_{2,3}^{k'}) \pi(z_{3,2}^{n'}) \xi' \rangle
            &\text{if }m=m'.
  \end{cases}
\end{IEEEeqnarray*}
Using the equalities $z_{3,3}^{*}=z_{1,1}^{}z_{2,2}^{}$,
$z_{2,2}^{}z_{2,3}^{}=0=z_{2,2}^{}z_{3,2}^{}$ and
$\pi(z_{2,2}^{})\xi=0=\pi(z_{2,2}^{})\xi'$, we get
\begin{IEEEeqnarray}{rCl}
  \langle \xi_{k,m,n}^{},  {\xi'}_{k',m',n'}^{} \rangle
  &=&
  \begin{cases}
    0 &\text{if }m\neq m',\\
    \langle  \pi(z_{2,3}^{k}) \pi(z_{3,2}^{n}) \xi \;, \;
    \pi(z_{2,3}^{k'}) \pi(z_{3,2}^{n'}) \xi' \rangle
        &\text{if }m=m'.
  \end{cases}
\end{IEEEeqnarray}
Next, note that from the equality
$z_{2,2}^{*}z_{2,2}^{}+z_{3,2}^{*} z_{3,2}^{} = 1$, it follows that
\[
(z_{3,2}^{*})^{m}(z_{3,2}^{n}) =
    (z_{3,2}^{*})^{m-1}(z_{3,2}^{n-1}) \quad
    \text{if }m\geq 1, n\geq 1,\text{ and } m+n>2.
\]
Therefore using the equalities  $\pi(z_{3,2}^{*})\xi=0=\pi(z_{3,2}^{*})\xi'$,
we get
\begin{IEEEeqnarray}{rCl}
  \langle \xi_{k,m,n}^{},  {\xi'}_{k',m',n'}^{} \rangle
  &=&
  \begin{cases}
    0 &\text{if }m\neq m' \text{ or }n\neq n',\\
    \langle  \pi(z_{2,3}^{k}) \xi \;, \;
    \pi(z_{2,3}^{k'}) \xi' \rangle
    &\text{if } m=m', n=n'.\\
  \end{cases}
\end{IEEEeqnarray}
Finally, note that
\[
z_{2,3}^{*}z_{2,3}^{}=z_{1,1}^{}z_{3,2}^{}z_{2,3}^{}=
z_{1,1}^{}z_{2,3}^{}z_{3,2}^{}=z_{3,2}^{*}z_{3,2}^{}.
\]
This together with the relations $z_{2,3}^{}z_{3,2}^{}=z_{3,2}^{}z_{2,3}^{}$
and $z_{2,3}^{}z_{3,2}^{*}=z_{3,2}^{*}z_{2,3}^{}$ give us
\begin{IEEEeqnarray}{rCl}
  \langle \xi_{k,m,n}^{},  {\xi'}_{k',m',n'}^{} \rangle
  &=&
  \begin{cases}
    0 &\text{if }m\neq m' \text{ or }n\neq n' \text{ or } k\neq k',\\
    \langle   \xi , \xi' \rangle &\text{if } m=m', n=n', k=k'.\\
  \end{cases}
\end{IEEEeqnarray}
Since $PQ=0$, we have $\langle \xi_{k,m,n}^{},  {\xi'}_{k',m',n'}^{} \rangle=0$.
\qed

The following two corollaries are now immediate.
\begin{corollary}\label{crlre:case6c}
There exists a $\lambda\in S^{1}$ such that the restriction of
$\pi(z_{1,3})$ to $P_{1,3}P_{3,1}\mathcal{H}$ is $\lambda I$.
\end{corollary}

\begin{corollary}\label{crlre:case6d}
There exists a $\mu\in S^{1}$ such that the restriction of
$\pi(z_{3,1})$ to $P_{1,3}P_{3,1}\mathcal{H}$ is $\mu I$.
\end{corollary}

\begin{proposition}\label{ppsn:case6e}
The representation $\pi$ is unitarily equivalent to
$\psi_{(\lambda,\bar{\mu}),s_{1}s_{2}s_{1}}$.
\end{proposition}
\prf
In view of the above corollaries, it follows from
(\ref{eq:action-6a-11}--\ref{eq:action-6a-33}) that
\begin{IEEEeqnarray}{rCl}
    \pi(z_{1,1}^{})\xi_{k,m,n}^{}
    &=&
    \begin{cases}
    \xi_{k-1,m,n-1}^{}&\text{if } k>0,\, n>0,\\
    0 & \text{otherwise,}\\
    \end{cases}\\
      &&\nonumber \\
     \pi(z_{1,2})\xi_{k,m,n}^{}
    &=&
     \begin{cases}
     \lambda\xi_{0,m-1,n+1}^{}&\text{if }k=0,\,m>0, \\
     \bar{\mu}\xi_{k-1,m,0}^{} &\text{if }k>0,\,n=0,\\
    0 & \text{otherwise,}
       \end{cases}\\
      &&\nonumber \\
   \pi(z_{1,3}^{}) \xi_{k,m,n}^{}
  &=&
  \begin{cases}
    \lambda\xi_{0,0,n}^{} &\text{if } m=0, \,k=0,\\
    0&\text{otherwise}, 
  \end{cases}\\
    &&\nonumber \\
   \pi(z_{2,1}^{})\xi_{k,m,n}^{}
  &=&
  \begin{cases}
    \bar{\lambda}\xi_{0,m,n-1}^{}&\text{if }k=0,\,n>0,\\
      \mu\xi_{k+1,m-1,0}^{}&\text{if }m>0,\,n=0,\\
    0 & \text{otherwise,}
    \end{cases}\\
    &&\nonumber \\
   \pi(z_{2,2}^{})  \xi_{k,m,n}^{}
  &=&
  \begin{cases}
    0 &\text{if  }m=0,\\
    \xi_{k+1,m-1,n+1}^{}  & \text{if  } m>0,
  \end{cases}\\
    &&\nonumber \\
   \pi(z_{2,3}^{})  \xi_{k,m,n}^{}
  &=&
  \begin{cases}
     \xi_{k+1,0,n}^{} &\text{if  }m=0,\\
     0 &\text{if  }m>0,
    \end{cases}\\
    &&\nonumber \\
   \pi(z_{3,1}^{})  \xi_{k,m,n}^{}
  &=&
  \begin{cases}
    \mu\xi_{k,0,0}^{} &\text{if  } m=0,\,n=0,\\
    0 &\text{otherwise,}
  \end{cases}\\
    &&\nonumber \\
   \pi(z_{3,2}^{})  \xi_{k,m,n}^{}
  &=&
  \begin{cases}
    \xi_{k,0,n+1}^{}  &\text{if  }m=0,\\
    0 &\text{if  }m>0,
  \end{cases}\\
    &&\nonumber \\
   \pi(z_{3,3}^{})  \xi_{k,m,n}^{} &=& \xi_{k,m+1,n}^{}.
\end{IEEEeqnarray}
Thus if we define $\mathcal{H}_{\xi}$ to be the clsoed linear span of
$\{\xi_{k,m,n}^{}:k,m,n\in\mathbb{N}\}$, then $\mathcal{H}_{\xi}$ is an
invariant subspace for $\pi$. Since $\pi$ is irreducible, we have
$\mathcal{H}=\mathcal{H}_{\xi}$. By the computations in the proof of the second
part of Proposition~\ref{ppsn:case6b}, it follows that
$\{\xi_{k,m,n}^{}:k,m,n\in\mathbb{N}\}$ is an orthonormal basis for
$\mathcal{H}$. The  map
\[
e_{m,k,n} \mapsto  \lambda^{k}\mu^{k-m-n} \xi_{k,m,n}^{}
\]
from $\ell^{2}(\mathbb{N}\times\mathbb{N}\times\mathbb{N})$ to $\mathcal{H}$ now
sets up a unitary equivalence between $\pi$ and
$\psi_{(\lambda,\bar{\mu}),s_{1}s_{2}s_{1}}$.
\qed

\section*{Acknowledgement}
The first author would like to thank the Indian Statistical Institute for supporting 
this work through a PhD fellowship.



\end{document}